\documentclass[leqno]{article}
\usepackage{a4wide,amssymb,latexsym,amsmath}
\newtheorem{itheorem}{Theorem}
\newtheorem{theorem}{Theorem}[section]

\newtheorem{prop}[theorem]{Proposition}
\newtheorem{corol}[theorem]{Corollary}
\newtheorem{defi}[theorem]{Definition}
\newtheorem{lemma}[theorem]{Lemma}

\def\noproof{{\unskip\nobreak\hfill\penalty50\hskip2em\hbox{}%
     \nobreak\hfill$\square$\parfillskip=0pt%
     \finalhyphendemerits=0\par}}
\def\enddemo{\ifmmode\eqno\square\else\noproof\vskip 12pt plus 3pt minus 9pt
\fi}
\def\diagram{\renewcommand\arraystretch{1.5} $$ \begin{array}}
\def\enddiagram{\end{array} $$ \renewcommand\arraystretch{1}}

\def\rmnewname#1{\expandafter\gdef\csname#1\endcsname{{\mathop{\rm
#1}\nolimits}}}
\def\itnewname#1{{\expandafter\gdef\csname#1\endcsname{{\mathop{\it
#1}\nolimits}}}}

 \itnewname{ab}
 \itnewname{ann}
 \itnewname{can}
 \rmnewname{cd}
 \itnewname{cl}
 \rmnewname{codim}
 \rmnewname{coker}
 \rmnewname{cone}
 \itnewname{cor}
 \rmnewname{CH}
 \rmnewname{diag}
 \rmnewname{div}
 \rmnewname{Div}
 \itnewname{et}
 \rmnewname{Et}
 \rmnewname{Ext}
 \itnewname{Frob}
 \itnewname{gen}
 \itnewname{gl}
 \rmnewname{GL}
 \itnewname{Gys}
 \itnewname{hor}
 \rmnewname{Hom}
 \itnewname{id}
 \rmnewname{im}
 \itnewname{incl}
 \rmnewname{Inj}
 \rmnewname{Iso}
 \rmnewname{Ker}
 \rmnewname{length}
 \itnewname{loc}
 \rmnewname{Mor}
 \rmnewname{NK}
 \rmnewname{ns}
 \rmnewname{Ob}
 \rmnewname{ord}
 \rmnewname{Pic}
 \rmnewname{Quot}
 \itnewname{rec}
 \rmnewname{red}
 \rmnewname{reg}
 \itnewname{res}
 \rmnewname{Sch}
 \itnewname{sh}
 \rmnewname{SK}
 \rmnewname{Sm}
 \rmnewname{sing}
 \rmnewname{SmCor}
 \itnewname{sp}
 \rmnewname{Spec}
 \rmnewname{Sub}
 \rmnewname{supp}
 \itnewname{tame}
 \rmnewname{Tame}
 \rmnewname{top}
 \rmnewname{Tor}
 \rmnewname{Tr}
 \itnewname{Zar}
 
 \def\C{{\mathbb C}}
 \def\Ce{{\mathcal C }}
 \def\E{{\mathcal E}}
 \def\F{{\mathbb F}}
 \def\G{{\mathbb G}}

  \def\M{{\mathcal M}}
 \def\N{{\mathbb N}}
 \def\O{{\mathcal O}}

 \def\p{{\mathfrak p}}

 \def\Q{{\mathbb Q}}
  \def\P{{\mathcal P}}
  
   \def\R{{\mathbb R}}
  \def\X{{\mathfrak X}}
  
  \def\D{{\mathfrak D}}
 \def\Z{{\mathbb Z}}

 \def\ms{\medskip}
 \def\ds{\displaystyle}

 \def\sst{\scriptstyle}
 \def\ssz{\scriptsize}

 \def\Cal#1{{\cal #1}}

 \def\lang{\longrightarrow}
 
 \def\mapr#1{\mathrel{\mbox{$\stackrel{#1}{\longrightarrow}$}}}
 \def\mapl#1{\mbox{$\stackrel{#1}{\longleftarrow}$}}
 \def\mapd#1{\Big\downarrow\rlap{$\vcenter{\hbox{$\scriptstyle{{#1}}$}}$}}
 
 \def\liso{\mathrel{\hbox{$\longrightarrow$} \kern-12pt\lower-4pt%
 \hbox{$\scriptstyle\sim$}\kern7pt}}
 \def\rightiso{\mathrel{\hbox{$\longleftarrow$} \kern-10pt\lower-4pt%
 \hbox{$\scriptstyle\sim$}\kern7pt}}
 \def\eqd{\Big\|}
 
 \def\mapu#1{\Big\uparrow\rlap {$\vcenter{\hbox{$\scriptstyle{{#1}}$}}$}}
 \def\mapd#1{\Big\downarrow\rlap {$\vcenter{\hbox{$\scriptstyle{{#1}}$}}$}}
 \def\surjr#1{\stackrel{#1}{\hbox{$\relbar \! \! \twoheadrightarrow$}}}
 
 \def\injr#1{\stackrel{#1}{\lhook\joinrel\relbar\joinrel\rightarrow}}
 
 \def\surjd#1{\lower4pt\hbox{$\downarrow$}\kern-.57em\Big\downarrow\rlap%
  {$\vcenter{\hbox{$\scriptstyle{{#1}}$}}$}}
 \def\surju#1{\lower-4pt\hbox{$\uparrow$}\kern-5.8pt\Big\uparrow\rlap%
  {$\vcenter{\hbox{$\scriptstyle{{#1}}$}}$}}

 \def\injd#1{
 \setlength{\unitlength}{0.1pt}
 \begin{picture}(40,120)
 \put(20,102){\vector(0,-1){158}} \put(13.5,100){\mbox{$\scriptscriptstyle
 \cap $}} \put(47,10){\mbox{$\scriptstyle #1 $}}
 \end{picture}}

 \def\inju#1{
 \setlength{\unitlength}{0.1pt}
 \begin{picture}(40,120)
 \put(20,-21){\vector(0,1){142}} \put(-16.3,-50){\mbox{$\scriptscriptstyle
 \cup $}} \put(47,10){\mbox{$\scriptstyle #1 $}}
 \end{picture}}

 \def\dirlim{\mathop{\lim\limits_{{\longrightarrow}}\,}}

\def\tcup{\mathop{\mathord{\cup}\mkern-8.5mu^{\hbox{.}}\mkern 9mu}}

 \def\hang{\hangindent\Itemindent}
 \def\textindent#1{\hskip\Itemindent\llap{\hbox{\rm #1}\enspace}\ignorespaces}
 \def\Item{\par\noindent\hang\textindent}
 
 \newdimen\Itemindent \Itemindent=.9cm

  \def\m{{\mathfrak m}}
 \def\Zt{{\mathcal Z}}
 \def\Ce{{\mathcal C}}

\date{20.9.2004}
\author{by Alexander Schmidt}
\title{\LARGE \bf  Tame class field theory for arithmetic schemes}

\begin{document}
\maketitle 
 The aim of global class field theory is the description of
abelian extensions of a finitely generated field $k$ in terms of its arithmetic invariants. The solution of this problem in the case of fields of dimension $1$ was one
of the major achievements of number theory in the first part of the previous century.

\medskip
In the 1980s, mainly due to K. Kato and S. Saito \cite{K-S3}, a generalization  to higher
dimensional fields has been found. The description of the abelian extensions is given in
terms of a generalized id\`ele class group, whose rather involved definition is based on
Milnor $K$-sheaves. However, if one restricts attention to unramified extensions (with
respect to a regular, projective model $X$ of $k$), then  class field theory has a nice
geometric description using the Chow group $\CH_0(X)$ of zero-cycles modulo rational
equivalence (see \cite{K-S1}, \cite{Sa}).

\bigskip
In the case of positive characteristic, a description of a similar geometric flavor of the
tamely ramified  abelian extensions of $k$ (w.r.t.\ a finite set $ D_1,\ldots, D_r$ of prime
divisors of $X$) was given in \cite{S-S}. The objective of this paper is to prove the
following result in the mixed characteristic case.

\begin{itheorem} \label{main}
Let $X$ be a regular connected scheme, flat and proper over $\Spec(\Z)$ such that its
generic fibre $X\otimes_\Z \Q$ is projective over $\Q$. Let $D$ be a divisor on
$X$ whose vertical irreducible components are normal schemes. Then there exists a natural isomorphism of finite abelian groups
\[
\rec_{X,D}: \CH_0(X,D) \liso \tilde \pi_1^t(X,D)^\ab\;.
\]
\end{itheorem}
Let us explain the ingredients of the theorem. First of all, $\tilde \pi_1^t(X,D)^\ab$ is the
abelianized modified tame fundamental group, which classifies finite abelian \'{e}tale coverings
of $X-\supp(D)$ which are at most tamely ramified along $D$ (cf.\ section~\ref{tamesect})
and in which every real point splits completely. $\CH_0(X,D)$ is the relative Chow group of
zero-cycles (cf.\ section~\ref{relk}). It is a quotient of the group of zero-cycles on $X-\supp(D)$
by an equivalence relation which is given by the divisors of functions satisfying a condition with respect to $D$. This equivalence relation is in general finer than rational equivalence.  Finally, the
map $\rec_{X,D}$ is uniquely determined by the property that the class of a closed point $x
\in X-\supp(D)$ is sent to its Frobenius automorphism in $\tilde \pi_1^t(X,D)^\ab$.

\medskip\noindent {\bf Remarks.} 
1. In many cases (e.g.\ if $X$ is semi-stable), the condition in the theorem on the vertical
components of $D$ is void. Furthermore, this condition can be weakened (see theorem~\ref{sharper}).

\smallskip\noindent
2. If $D$ is zero, Theorem~\ref{main} reduces to the ``unramified class field theory'' for
arithmetic schemes of K. Kato and S. Saito \cite{K-S1}, \cite{Sa} (which we use in the
proof).

\smallskip\noindent
3. The finiteness of the group $\tilde \pi_1^t(X,D)^\ab$ was shown in \cite{S1}. Furthermore
(loc.cit.) this group depends only on the scheme $U=X-\supp(D)$. Hence it is desirable to
give a formulation of the tame class field theory of $U$ solely in terms of $U$. We
conjecture that the above theorem is true with the relative Chow group replaced by $h_0(U)$,
the 0th singular homology group of $U$, which was considered in \cite{S2}. Indeed, there
exists a surjective reciprocity map and the  conjecture is that this surjection is an
isomorphism. At the moment this is known if $\dim X=1$ or if $D$ is zero.

\smallskip\noindent
4. If $\dim X=1$, then $\CH_0(X,D)$ is isomorphic to the ray class group of the number field
$k(X)$ with modulus $\m_D$ where $\m_D$ is the (square-free) product of the points in $D$.
In this case, $\rec$ is the reciprocity homomorphism of  classical (one-dimensional)
class field theory.

\vskip.8cm I want to thank M.\ Spie{\ss} for many remarks and comments he made during our joint
work \cite{S-S}, which was the starting point of present paper. I am grateful to U.\ Jannsen
for several helpful discussions and in particular for calling my attention to the existence
of \cite{Po}. Hearty thanks go to M.\ Levine from whom I learned much about relative
$K$-groups.  Finally, I  want to thank D.\ Grayson for his comments on a preliminary version
of this paper. A large parts of this work was written while I was a Heisenberg fellow of the Deutsche
Forschungsgemeinschaft and I want to thank this institution for their financial support.

\vskip.8cm\noindent {\it{\bf Notational Convention:} Throughout this paper,  given a scheme\/
$Y$, we denote by $Y^+$ the disjoint union of its irreducible components and by $\tilde{Y}$
the normalization in its total ring of fractions.}

\section{Relative \boldmath Chow groups} \label{relk}

We start by collecting some facts on relative $K$-groups. A more detailed discussion and
proofs of the results mentioned below are found in \cite{Le}, \S1. Following Quillen
\cite{Qu}, the $K$-groups of an exact category $\Ce$ are defined by
\[
K_i\Ce=\pi_{i+1} (K\Ce,*), \quad K\Ce= BQ\Ce.
\]
Here $B$ means geometric realization, $Q$ is Quillen's construction \cite{Qu} and
$\pi_i(-,*)$ is the $i$-th homotopy group of a pointed topological space.

\medskip
For a noetherian scheme $X$ let $\M_X$ be the category of coherent $\O_X$-sheaves and
$G_iX=K_i\M_X$. Assume that $X$ is regular and let $D=D_1+\cdots+ D_r$ be a sum of
distinct irreducible divisors on $X$.  Let $Y:=\supp(D)=D_1 \cup \cdots \cup D_r \subset X$, hence
$Y^+= D_1 \tcup \cdots \tcup D_r$. Let $\M_{X|Y}$ be the full subcategory in $\M_X$
consisting of coherent $\O_X$-sheaves $F$ with $\Tor^{\O_X}_i(F,\O_Y)=0$ for $i>0$. Since
every coherent $\O_X$-sheaf is the quotient of a vector bundle, the resolution theorem shows
that the inclusion $\M_{X|Y} \to \M_X$ induces an isomorphism on $K$-groups. The natural
functor $j_{Y^+}^*: \M_{X|Y} \to \M_{Y^+}$ is exact and we define the {\boldmath \bf
relative $G$-groups $G_*(X,D)$} as the homotopy groups of the homotopy fibre of $j_{Y^+}^*:
K\M_{X|Y}\to K\M_{Y^+}$. If $f:X' \to X$ is a flat morphism inducing a flat morphism
$f_{|{Y'}^+}: {Y'}^+ \to Y^+$, then we obtain a natural homomorphism
\[
i^*: G_*(X,D) \lang G_*(X',D')
\]
(see \cite{Co}, \cite{Le} for the technical details). Recall  \cite{Qu} the usual Quillen
spectral sequence for $X$
\[
E_1^{pq}(X)= \bigoplus_{x\in X^p} K_{-p-q}(k(x))\Rightarrow G_{-p-q}(X). \leqno{(1)}
\]
which is associated to the filtration by codimension of support on $\M_X$
\[
\M_X^i=\{ F \in \M_X\,|\, \codim_X \supp (F) \geq i\}.
\]
As is well known, the Chow group of codimension $p$ cycles on $X$ and the term $E_2^{p,-p}$
of the above spectral sequence are naturally isomorphic.  Our next aim is to construct a
natural spectral sequence converging to $G_*(X,D)$. The relative Chow group of codimension
$p$ cycles will be the term $E_2^{p,-p}$ of that spectral sequence by definition.

Let $\M^i_{X|Y} \subset \M_{X|Y}$ be the full subcategory of coherent sheaves $F$ in
$\M_{X|Y}$ with\break $\codim_X\supp(F)\geq i$. Let $G^i(X,D)$ be the homotopy fibre of
$BQj_{Y^+}^*: BQ\M_{X|Y}^i \to BQ\M_{Y^+}^i$ and $G_p^i(X,D)=\pi_{p+1}(G^i(X,D))$. For $i
<k$, let $\M_{X|Y}^{i/k}$ be the direct limit
\[
\M_{X|Y}^{i/k}=
\dirlim_{\sst \renewcommand{\arraystretch}{0.8}\begin{array}{c} \sst Z \subset X\\
\hbox{\ssz $Z$ intersecting}\\
\hbox{\ssz $Y$ properly}\\ \sst \codim_XZ \geq k
\end{array}
\renewcommand{\arraystretch}{1}}
\M_{X-Z|Y-Z}^i.
\]
and let $\M_{Y^+}^{i/k}$ be the direct limit
\[
\M_{Y^+}^{i/k}= \dirlim_{\sst \renewcommand{\arraystretch}{0.8}\begin{array}{c}
\sst Z \subset Y^+\\
\sst \codim_{Y^+} Z \geq k
\end{array}
\renewcommand{\arraystretch}{1}}
\M_{Y^+-Z}^i.
\]
Let $G^{i/k}(X,D)$ be the homotopy fibre of $BQ\M^{i/k}_{X|Y}\to BQ\M^{i/k}_{Y^+}$ and put
$G_p^{i/k}(X,D) =\break\pi_{p+1}(G^{i/k}(X,D))$.  By \cite{Le}, 1.6, we obtain a
spectral sequence
\[
E_1^{pq}(X|Y) \Longrightarrow G_{-p-q}(X),\leqno{(2)}
\]
in which the $E_1$-terms are given by
\[
E_1^{pq}(X|Y)=\left\{
\begin{array}{ll}
K_{-p-q}(\M^{p/p+1}_{X|Y}),& -p-q>0,\ p\leq \dim X\\
\bar K_0(\M^{p/p+1}_{X|Y}), &-p-q=0\\
\; 0 & \hbox{ otherwise.}
\end{array}
\right.
\]
Here $\bar K_0(\M^{p/p+1}_{X|Y}):= \im \left( K_0(\M^{p}_{X|Y}) \lang
K_0(\M^{p/p+1}_{X|Y})\right)$. The natural functor $\M_{X|Y}^{p/p+1} \to
\M_X^{p/p+1}$ induces a map from the $E_1$-terms of this spectral sequence (2) to that of the
usual Quillen spectral sequence (1) which is compatible with differentials.  Furthermore,
one obtains a spectral sequence
\[
E_1^{pq}(X,D) \Longrightarrow G_{-p-q}(X,D), \leqno{(3)}
\]
in which the $E_1$-terms are given by
\[
\renewcommand{\arraystretch}{1.3    }
 E_1^{pq}(X,D)=\left\{
\begin{array}{ll}
G_{-p-q}^{p/p+1}(X,D),& -p-q>0,\ p\leq \dim X\\
\bar G_0^{p/p+1}(X,D), &-p-q=0\\
\; 0 & \hbox{ otherwise.}
\end{array}
\right.
\renewcommand{\arraystretch}{1}
\]
Here $\bar G_0^{p/p+1}(X,D)= \im \left(G_0^p(X,D) \to G_0^{p/p+1}(X,D)\right)$. The
filtration on $G_*(X,D)$ is the ``topological'' filtration:
\[
F^pG_*(X,D)= \im \left(  G_*^p(X,D) \lang G_*(X,D)\right).
\]

\begin{defi}
We call the group
\[
\CH^p(X,D)=E_2^{p,-p}(X,D)
\]
the relative Chow group of codimension $p$-cycles. If\/ $X$ is equidimensional of dimension
$d$, we call the group $ \CH_p(X,D)=E_2^{d-p,-d+p}(X,D)$ the relative Chow group of
$p$-cycles.
\end{defi}

\noindent
{\bf Remark:}  If \/ $D$ is empty, then
\[
\CH^p(X,\varnothing) = \coker\Big(\bigoplus_{x\in X^{p-1}} k(x)^\times \to \bigoplus_{x\in
X^p} \Z\Big)
\]
is the usual group of codimension $p$ cycles on $X$ modulo rational equivalence.

\bigskip Finally, we mention the existence of a long exact sequence of $E_1$-terms:
\[
\cdots \to E_1^{p,-p-2}(Y^+) \to E_1^{p,-p-1}(X,D) \to E_1^{p,-p-1}(X|Y) \to
E_1^{p,-p-1}(Y^+)
\]
compatible with the differentials. If the maps $G_0^p(X,D)\to G_0^{p/p+1}(X,D)$ and
$K_0(\M^{p}_{X|Y}) \to K_0(\M^{p/p+1}_{X|Y})$ are surjective, we can extend this
sequence to the right up to $E_1^{p,-p}(Y^+)$.

\section{Explicit description of \boldmath $\CH_0(X,D)$} \label{explsect}

The rather implicit definition of the relative Chow groups in the last section was given in
order to obtain good functorial properties. In this section we will give a more explicit
description of the group $\CH_0(X,D)$ as the group of zero cycles on $X-\supp(D)$ modulo
relations. These relation come via the divisor map from families of rational functions on curves in $X-\supp(D)$ which satisfy an additional condition. Informally speaking, the condition is that the product of the values of these functions at  points of $D$ must be $1$.

\medskip
 We assume from now on that $X$ is an excellent regular noetherian scheme. Then the
closed subscheme $Y=\supp(D)$ is locally principal,  locally defined by a non-zero divisor.
A coherent sheaf $F$ is in $\M_{X|Y}$ if and only if it is without $I_Y$-torsion and for
$F\in \M_{X|Y}$ the support $\supp(F)$ has proper intersection with $Y$.  If $Z$ is a closed
reduced subscheme, then $\O_Z$ is in $\M_{X|Y}$ if and only if $Z$ intersects $Y$ properly.
For $0 \leq p \leq \dim (X)$, we consider the (abelian) categories
\medskip
\Itemindent=2cm
\Item{$\M_{X:Y}^p -$}
coherent $\O_X$-sheaves whose support is contained in a subscheme of codimension $\geq p$
which has proper intersection with $Y$.

\medskip\noindent
We have a natural inclusion $\M_{X|Y}^p \to \M_{X:Y}^p$.

\begin{lemma}\label{zweitess}
The natural map $Q\M_{X|Y}^p \to Q\M_{X:Y}^p$ is a homotopy equivalence for all $p$.
\end{lemma}

\begin{demo}{Proof:}
Let $F\in \M_{X|Y}^p$ and let $Z$ be the (not necessarily reduced) subscheme associated with
the ideal sheaf $\ann(F)$. Then $F$ is an $\O_Z$-module and $\O_Z \in \M_{X|Y}^p$.
Therefore the natural inclusion $\M_{X|Y}^p \to \M_{X:Y}^p$ can be factored as
\[
\M_{X|Y}^p \mapr{i} C^p \mapr{j} \M_{X:Y}^p,
\]
where $C^p$ is the full subcategory of objects in $\M_{X:Y}^p$ which are $\O_Z$-modules for
some closed (not necessarily reduced) subscheme $Z\subset X$ with $\O_Z \in \M_{X|Y}^p$. The
category $C^p$ is closed under taking subobjects, quotient and finite direct products and is
therefore abelian.

Let $F$ be a coherent $\O_X$-sheaf whose support is of codimension $\geq p$ and has proper
intersection with $Y$. Since $X$ is noetherian, the coherent ideal sheaf $\ann(F)$ contains a
power of its radical. Therefore we find a filtration $0=F_0 \subset F_1 \subset \cdots
\subset F_n=F$ of $F$ such that all subquotients are $\O_Z$-modules, where $Z=\supp(F)$ as
reduced subscheme. Hence all subquotients are in $C^p$ and we conclude by devissage that
$Qj$ is a homotopy equivalence.

Now assume that $F \in C^p$, i.e.\ $F$ is an $\O_Z$-module where $Z$ is a closed subscheme
with $\O_Z \in \M_{X|Y}^p$. Since $X$ is regular, we find a locally free $\O_X$-sheaf $\cal
F$ and a surjection ${\cal F} \twoheadrightarrow F$. It factors to a surjection
\[
 {\cal F}/I_Z \surjr{p} F.
\]
The sheaf ${\cal F}/I_Z$ (being a locally free $\O_Z$-sheaf) is without $I_Y$-torsion.
Therefore ${\cal F}/I_Z$ and also its subsheaf $\ker(p)$ are in $\M_{X|Y}^p$ and Quillen's
resolution lemma implies that $Qi$ is a homotopy equivalence.
\end{demo}

For a closed subscheme $Z\subset X$, let $\M_{X}(Z)$ be the category of coherent
$\O_X$-sheaves with support in $Z$. Defining $\M_{X:Y}^{p/p+1}$ as the direct limit
\[
\M_{X:Y}^{p/p+1}= \dirlim_{\sst \renewcommand{\arraystretch}{0.8}
\begin{array}{c} \sst W \subset Z \subset X\\
 \sst \codim_XZ \geq p\\
\sst \codim_X W \geq p+1\\
\hbox{\ssz $W$,$Z$  intersect $Y$ properly}
\end{array}
\renewcommand{\arraystretch}{1}}
\M_{X-W}(Z),
\]
we conclude that also the natural map $Q\M_{X|Y}^{p/p+1} \to Q\M_{X:Y}^{p/p+1}$ is a
homotopy equivalence. This allows to calculate some terms of the spectral sequence (2) of
the last section.

\begin{prop}\label{ausrech}
The map $i_{X-Y}^*: \M_{X|Y}^{p/p+1} \to \M_{X-Y}^{p/p+1}$ induces an isomorphism
\[
K_0(\M_{X|Y}^{p/p+1}) \liso K_0(\M_{X-Y}^{p/p+1}) \liso \bigoplus_{x \in (X-Y)^p} \Z \; .
\]
For $q= 1,2$, we obtain short exact sequences

\diagram{ccccccccc}
 0 &\to& K_q(\M_{X|Y}^{p/p+1})& \to& K_q(\M_{X-Y}^{p/p+1})& \stackrel{\delta}{\to}&
K_{q-1}(\M_{Y}^{p/p+1}) &\to &0\\
&&&&\mapu{\wr}&&\mapu{\wr} &&\\
&&&&\ds\bigoplus_{x \in (X-Y)^p}\!\!\!\!K_q(k(x))&&\ds\bigoplus_{y \in Y^p} K_{q-1}(k(y)).
\enddiagram
Here $\delta$ is the composition
\diagram{ccccccc}
K_q(\M_{X-Y}^{p/p+1})&\hspace{-1cm}\mapr{\incl}\hspace{-.6cm}&
K_q(\M_{X}^{p/p+1})&\hspace{-.5cm}\mapr{\partial}\hspace{-.5cm}
&K_{q-1}(\M_{X}^{p+1/p+2})&\hspace{-.5cm}\mapr{proj}\hspace{-.5cm}&K_{q-1}(\M_{Y}^{p/p+1})\\
\mapu{\wr}&&\mapu{\wr}&&\mapu{\wr}&&\mapu{\wr}\\
\ds\bigoplus_{x \in (X-Y)^p}\!\!\!\!K_q(k(x))&&\ds\bigoplus_{x \in
X^p}K_q(k(x))&&\ds\bigoplus_{x \in X^{p+1}} K_{q-1}(k(x))&&\ds\bigoplus_{y \in
Y^p} K_{q-1}(k(y)),
\enddiagram
where $\partial$ is the boundary map of the Quillen spectral sequence.
\end{prop}

\begin{demo}{Proof:} First of all, we can replace $\M_{X|Y}^{p/p+1}$ by $\M_{X:Y}^{p/p+1}$.
The localization theorem shows that
\[
Q\M_{Y}^{p/p+1} \to Q\M_{X:Y}^{p/p+1} \to Q\M_{X-Y}^{p/p+1}
\]
is a homotopy fibre sequence. This implies a long exact sequence
\[
 \cdots \to K_1 (\M_{X-Y}^{p/p+1}) \stackrel{\delta_1}{\to} K_0(\M_{Y}^{p/p+1})
\to K_0(\M_{X:Y}^{p/p+1}) \to K_0(\M_{X-Y}^{p/p+1})
\]
Let $x \in (X-Y)^p$ and let  $Z=\overline{ \{x\}} $ be its closure in $X$. Then $Z$ has
proper intersection with $Y$ and the class of the free rank one  $\O_Z$-module in
$K_0(\M^p_{X:Y})$ maps to the class of the one-dimensional $k(x)$-vector space in
$K_0(\M_{X-Y}^{p/p+1})$. Hence $K_0(\M_{X:Y}^{p/p+1}) \to K_0(\M_{X-Y}^{p/p+1})$ is
surjective and in order to conclude the proof, it remains to show that the maps
$\delta_q:K_q(\M_{X-Y}^{p/p+1}) \to K_{q-1}(\M_{Y}^{p/p+1})$ are surjective for $q \leq 3$.

As is well known, $K_{q-1}(\M_Y^{p/p+1})=\bigoplus_{y \in Y^p} K_{q-1}(k(y))$.  Let $y \in
Y^p$ be any point.  It suffices to show that any element of $\alpha\in K_{q-1}(k(y))$
(considered as an element of $K_{q-1}(\M_Y^{p/p+1})$) is in the image of $\delta_q$.
Furthermore ($q\leq 3$), we may suppose that $\alpha=\{a_1,\ldots, a_{q-1}\}$ is a symbol.
Choose a closed integral subscheme $Z$ of codimension $p$ in $X$ which has proper
intersection with $Y$ and such that $y$ is a regular point of $Z$. Let $R$ be the semi-local
ring of $Y\cap Z$ in $Z$ and let $F$ be its quotient field. Since $X$ is excellent, the
normalization $\tilde R$ of $R$ in $F$ (a semi-local PID) is finite over $R$. So we obtain
the following commutative diagram

\[
\renewcommand{\arraystretch}{0.8}
\begin{array}{ccc}
K_q(\M_{X-Y}^{p/p+1})& \lang& K_{q-1}(\M_Y^{p+1})\\
\mapu{}&&\mapu{}\\
K_q(F)&\lang& \ds \bigoplus_{\m \subset R} K_{q-1}(R/\m)\\
\mapu{\wr}&&\mapu{}\\
K_q(F)&\lang& \ds \bigoplus_{\tilde\m \subset \tilde R} K_{q-1}(\tilde R/\tilde \m)
\end{array}
\renewcommand{\arraystretch}{1}
\]
Since $y$ is a regular point of $Z$, there is exactly one $\tilde \m_0 \subset \tilde R$
corresponding to $y$ and $\tilde R/\tilde m_0 =k(y)$. Now let $\pi \in \tilde R$ be an
element which is a local parameter for $\tilde \m_0$ and is congruent to $1$ modulo all other
$\tilde \m \subset \tilde R$. Then the symbol $\{a_1,\ldots, a_{q-1},\pi\} \in K_q(F)$ defines
the required pre-image of $\alpha$ in $K_q(\M_{X-Y}^{p/p+1})$.
\end{demo}

\begin{lemma} \label{3.6} If $X$ is pure of dimension $d$, then
\[
E_1^{d,-d}(X,D) \cong \bigoplus_{x \in (X-Y)^d} \Z\; .
\]
\end{lemma}

\begin{demo}{Proof:}
By definition, $E_1^{d,-d}(X,D)= \im(G_0^d(X,D) \to G_0^{d/d+1}(X,D))$. By dimension reasons
and by proposition~\ref{ausrech}, we have a commutative diagram
\[
\begin{array}{ccccc}
G_0^d(X,D)& \liso & K_0(\M_{X|Y}^{d})\\
\mapd{}&&\mapd{\alpha}&\\ G_0^{d/d+1}(X,D)& \liso &
K_0(\M_{X|Y}^{d/d+1})&\liso&\ds\bigoplus_{x \in (X-Y)^d} \Z\;.
\end{array}
\]
Furthermore, $\alpha$ is surjective, since it has  a section $s_\alpha$ which is defined by
sending the class of an $n$-dimensional vector space over $k(x)$, $x\in (X-Y)^p$ to the
class of the free rank-$n$ module over $\O_{\bar{\{x\}}}$. This implies the statement.
\end{demo}

\begin{prop}\label{chow0descrklein} Assume that $X$ is a noetherian, excellent, regular
and connected scheme of pure dimension $d$.  Let $D$ be a divisor on $X$ and let
$Y=\supp(D)$. Then the group $\CH_0(X,D)$ is the quotient of $\ds\bigoplus_{x\in (X-Y)^d}\Z$
by the image of the group
\[
R_{X,D}\stackrel{\rm df}{=}\ker \Big( \ker \big(\bigoplus_{x \in (X-Y)^{d-1}} k(x)^\times \to
\bigoplus_{y \in Y^{d-1}} \Z \big) \lang  \bigoplus_{z \in {(Y^+)}^{d-1}} k(z)^\times \Big)
\]
under the divisor map.
\end{prop}

\begin{demo}{Proof:} By definition,
$
\CH_0(X,D)=E_2^{d,-d}(X,D)= E_1^{d,-d}(X,D)/\im\, E_1^{d-1,d}(X,D)
$.
By lemma~\ref{3.6}, we have a commutative diagram with exact rows
\[
\renewcommand{\arraystretch}{1.6}
\begin{array}{ccccccc}
0& \to & E_1^{d,-d}(X,D)& \stackrel{\sim}{\to}& E_1^{d,-d}(X|Y)& \to& 0\\
\mapu{}&& \mapu{}&& \mapu{}&& \mapu{} \\
 E_1^{d-1,-d-1}(Y^+) &\to& E_1^{d-1,-d}(X,D)& \to& E_1^{d-1,-d}(X|Y)& \to& E_1^{d-1,-d}(Y^+)
\end{array}
\renewcommand{\arraystretch}{1}
\]
Therefore, $\CH_0(X,D)$ is the quotient of $E_1^{d,-d}(X,D)\cong\bigoplus_{x\in (X-Y)^d}\Z $
by the image of\break $\ker(E_1^{d-1,-d}(X|Y) \to E_1^{d-1,-d}(Y^+))$. The statement of the
proposition follows from proposition~\ref{ausrech}.
\end{demo}

In order to understand the group $\CH_0(X,D)$ explicitly, it remains to give a description of
the maps
\[
\delta_0: \bigoplus_{x \in (X-Y)^{d-1}} k(x)^\times \to \bigoplus_{y \in Y^{d-1}} \Z
\]
and
\[
\phi: \ker \big(\bigoplus_{x \in (X-Y)^{d-1}} k(x)^\times \stackrel{\delta_0}{\to}
\bigoplus_{y \in Y^{d-1}} \Z \big) \lang  \bigoplus_{z \in {(Y^+)}^{d-1}} k(z)^\times
\]
of the last proposition. The map $\delta_0$ is determined by the boundary map $\partial_0:
K_1(\M_{X}^{d-1/d}) \to K_0(\M_X^{d/d+1})$ of the Quillen spectral sequence and an explicit
description of $\partial_0$ is given at the beginning of the appendix. To describe $\phi$,
we need the following lemma which describes the natural map $\res: E_1^{p,-p-q}(X|Y) \to
E_1^{p,-p-q}(Y^+)$ if $X$ is affine and $Y=D$ is irreducible.

\begin{lemma} \label{commdiag}
If\/ $X=\Spec(A)$ is affine and $Y$ is irreducible, defined by a non-zero element $r \in A$,
then the following diagram commutes for all $p,q\geq 0$
\diagram{ccc} K_q(\M_{X|Y}^{p/p+1})&\mapr{\can} &K_q(\M_{X-Y}^{p/p+1})\\
\mapd{\res}&&\mapd{\cdot \cup r}\\
K_q(\M_Y^{p/p+1})& \mapl{\delta_q}& K_{q+1}(\M_{X-Y}^{p/p+1}).
\enddiagram
Here $\res$ is the natural map induced by the embedding of\/ $Y$ to $X$, $\delta_q$ is the
composite map $K_{q+1}(\M_{X-Y}^{p/p+1}) \hookrightarrow K_{q+1}(\M_{X}^{p/p+1})
\stackrel{\partial_q}{\to} K_{q}(\M_{X}^{p+1/p+2}) \stackrel{pr}{\to} K_q(\M_Y^{p/p+1})$ and
$\cdot \cup r$ is the cup product with the class of $r$ in $K_1(A[r^{-1}])=K_1(X-Y)$, induced
by the pairing
\[
\P(X-Y)\times \M_{X-Y}^{p/p+1}\to \M_{X-Y}^{p/p+1}, \quad (P,M)\mapsto P \otimes M.
\]
\end{lemma}

To describe $\phi$, it suffices to understand its
$z$-component for each $z \in Y^+$. Therefore we may suppose that $X=\Spec(A)$ is affine,
local and that $D=Y$ is irreducible, given as the zero locus of an element $0\neq r \in A$.
Let $y$ be the closed point of $Y$. By the lemma (for $p=d-1$ and q=1), we obtain the
following explicit description of the map
$$\phi: \ker \big(\bigoplus_{x \in (X-Y)^{d-1}} k(x)^\times \stackrel{\delta_0}{\to}
\Z \big) \lang  k(y)^\times.$$ An element $a=(a_x)_x\in \ker \big(\bigoplus_{x \in
(X-Y)^{d-1}} k(x)^\times \stackrel{\delta_0 }{\to} \Z\big)$ maps to the image of the symbol
$\{a,r\}=(\{a_x,r\})_x \in \bigoplus_{x \in (X-Y)^{d-1}} K_2(k(x))$ under the boundary map
$$\delta_1: \bigoplus_{x \in (X-Y)^{d-1}} K_2(k(x)) \lang k(y)^\times.$$ An explicit
description of the boundary map $\delta_1$ is given at the beginning of the appendix.

\begin{demo}{Proof of lemma~\ref{commdiag}:}
Let us fix some notation: For a category $\Ce$ let
its subdivision $\Sub(\Ce)$ be the category whose objects are morphisms in $\Ce$ and a
morphism from $f$ to $g$ is a commutative diagram
\[
\begin{array}{cccc}
\bullet &\mapr{f}&\bullet\\
\mapu{}&&\mapd{}\\
\bullet&\mapr{g}&\bullet&.
\end{array}
\]
The functor $\Sub(\Ce)\to \Ce$ which sends an arrow to its target is a homotopy equivalence
(\cite{Gr}, p.228). Assume that $\Ce$ is an exact category in which every exact sequence
splits. Let $\Inj(\Ce)$ and $\Iso(\Ce)$ be the categories whose objects are those of $\Ce$
and whose morphisms are the admissible monomorphisms and the isomorphisms in $\Ce$,
respectively. Put $I=\Iso(\Ce)$ and $\E=\Sub(\Inj(\Ce))$. Sending $M
\stackrel{i}{\hookrightarrow} N$ to $\coker(i)$ defines a natural functor $\E \to Q\Ce$. The
symmetric monoidal category $I$ acts on $\E$ by $P\oplus(M \stackrel{i}{\hookrightarrow} N)=
(P\oplus M \stackrel{\id_P\oplus i }{\hookrightarrow} P\oplus N)$ and we obtain a fibration
sequence (\cite{Gr}, p.228)
\[
I^{-1}I \to I^{-1}\E \to Q\Ce
\]
Furthermore, $I^{-1}\E$ is contractible which induces a homotopy equivalence
\[
\psi: BI^{-1}I \liso \Omega B Q\Ce
\]
(the well-known comparison map between Ext- and Q-construction).  For $M \in \Ob(\Ce)$, we
have two natural maps $(i_M)_!, (p_M)^!: 0 \to M$ in $Q\Ce$ which are associated to the
canonical maps $i_M: 0 \hookrightarrow M$ and $p_M: M \twoheadrightarrow 0$ in $\Ce$. The
map $\psi$ is characterized by sending an object $(N,M)$ of $I^{-1}I$ to the loop
\[
0 \mapr{(p_N)_!} N \mapl{(i_N)^!} 0 \mapr{(i_M)^!} M \mapl{(P_M)_!} 0
\]
in $BQ\Ce$.

Next, following \cite{Sn}, we construct a categorical description of the boundary map. Let
$\Ce=\M_{X|Y}^{p/p+1}$ (note that every exact sequence in $\Ce$ splits). Consider the
category with the same objects as $\Ce$ and whose morphisms are (not necessarily admissible) injections
$i: M \hookrightarrow N$ in $\Ce$ whose cokernel is annihilated by $r$, and let $H$ be its
subdivision. The category $H$ is the source of two functors:

\medskip\noindent
1. Let $J=\Iso(\M_{X-Y}^{p/p+1})$. Sending the object $i:M \hookrightarrow N$ of $H$ to the
restriction of $N$ to $X-Y$ yields a natural functor $ H \to J$.\\
2.  Sending $i:M \hookrightarrow N$ to $\coker(i)$ yields a functor $H\to Q\M_Y^{p/p+1}$.

\medskip\noindent
There is a natural $I$-action on $H$ given by $P\oplus(M \stackrel{i}{\hookrightarrow} N)=
(P\oplus M \stackrel{\id_P\oplus i}{\hookrightarrow} P\oplus N)$. We therefore obtain
functors
\[
\alpha: I^{-1}H \to J^{-1}J,\quad \beta: I^{-1}H \to Q\M_Y^{p/p+1}.
\]
By \cite{Sn}, proof of theorem 2.2 (our $H$ corresponds to $g^{-1}(0)$ there), $I^{-1}H \to
I^{-1}J$ is a homotopy equivalence and $I^{-1}J \to J^{-1}J$ is a homotopy equivalence on
base point components. Furthermore (loc.cit.),  the following diagram commutes for $q\geq 0$
\diagram{ccc}
\pi_{q+1} (BI^{-1}H )&\mapr{\beta_*}& \pi_{q+1}(B Q\M_Y^{p/p+1})\\
\mapd{\wr\:\alpha_*} && \eqd\\
\pi_{q+1} (B J^{-1}J)& \mapr{\delta}& \pi_{q+1}(B Q\M_Y^{p/p+1}),
\enddiagram
which gives a categorical description of the boundary map $\delta:
K_{q+1}(\M_{X-Y}^{p/p+1})\to K_q(\M_Y^{p/p+1})$.

\medskip
For an object $M$ of $I$ we have two morphism $r_M,r^M: (0,0=0) \to (0,r M \subset M) $ in
$I^{-1}H$ given by the commutative diagrams
\[
\begin{array}{cccl}
M&\mapr{\id_M}&M\\
\inju{}&&\mapd{\id_M}\\
r M &\injr{} &M&
\end{array} \qquad
\hspace*{.6cm} \hbox{\rm and}\hspace*{1cm}
\begin{array}{cccl}
M&\mapr{\id_M}&M\\
\mapu{(\cdot r)^{-1}}&&\mapd{\cdot r}\\
 r M &\injr{} &M&,
\end{array}
\]
where $(\cdot r)^{-1}$ is the inverse of the isomorphism $M\stackrel{\cdot r}{\to} rM$.
Consider the map $\gamma:BI^{-1}I \to \Omega BI^{-1}H$, which is characterized by sending an
object $(N,M)$ of $I^{-1}I$ to the loop
\[
(0,0=0) \mapr{r^N} (0,r N \subset N)\mapl{r_N} (0,0=0) \mapr{r_M} (0,r M \subset M)
\mapl{r^M} (0,0=0)
\]
in $BI^{-1}H$. Considering the natural image of  $r \in A[r^{-1}]^\times$ in the group
$K_1(A[r^{-1}])=\break\pi_1(B\GL(A[r^{-1}])^+)$ and using the comparison between, $+$ and
Ext-construction as described in \cite{Gr}, p.224,
we see that for $q \geq 0$ the induced map $\pi_{q}(BI^{-1}I) \mapr{\gamma} \pi_{q+1}(
BI^{-1}H) \liso \pi_{q+1}(BJ^{-1}J)$ is the product map
\[
\cdot \cup r: K_q(\M_{X|Y}^{p/p+1}) \lang K_{q+1}(\M_{X-Y}^{p/p+1}).
\]
Furthermore, $\Omega\beta \circ \gamma=\Omega\res\circ\psi$ where $\res: Q\M_{X|Y}^{p/p+1}
\to Q\M_Y^{p/p+1}$ is the natural restriction map. Summing up, we obtain a commutative
diagram
\diagram{ccc}
\Omega B I^{-1}H&\mapl{\gamma}&BI^{-1}I\\
\mapd{\Omega\alpha}&&\mapd{\Omega \res\circ \psi}\\
\Omega BJ^{-1}J& \mapr{\delta}& \Omega BQ\M_Y^{p/p+1}.
\enddiagram
Passing to homotopy groups, we obtain the statement of lemma~\ref{commdiag}.
\end{demo}

\section{Finiteness of \boldmath $\CH_0(X,D)$} \label{finsect}
The aim of this section is to detect relations in $\CH_0(X,D)$ and to use these relations to
show that $\CH_0(X,D)$ is finite, when $X$ is an arithmetic scheme. We first recall some well-known facts on the filtration on the $K$-groups $K_r(F)$ for $r=1,2$ of discrete valuation fields. Let $(F,v)$ be a discrete valuation field with residue field $k$. Let $\O_F$ be the valuation ring in $F$ and let $\m
\subset \O_F$ be its maximal ideal, thus $\O_F/\m \cong k$. Then one puts
$U^0(F^\times)=\O_F^\times$ and for $i \geq 1$, $U^i(F^\times)=\ker(\O_F^\times \to
(\O_F/\m^i)^\times)$. As is well known, we have natural isomorphisms
\[
F^\times/U^0(F^\times) \cong \Z\;, \qquad U^0(F^\times)/U^1(F^\times)\cong k^\times.
\]
For $i\geq 1$ the group $U^iK_2(F)$ is the subgroup generated by symbols $\{u,x\}$ with $u\in
U^i(F^\times)$ and $x \in F^\times$. The group $U^0K_2(F)$ is the kernel of the tame symbol
$K_2(F) \to k^\times$. We have inclusions
\[
K_2(F) \supseteq U^0K_2(F) \supseteq U^1K_2(F) \supseteq \cdots
\]
and natural isomorphisms (\cite{B-T}, prop.4.3,4.5)
\[
K_2(F)/U^0K_2(F)\cong k^\times\; , \qquad U^0K_2(F)/U^1K_2(F) \cong K_2(k),
\]
If $F_v$ and $\hat F_v$ denote the henselization and the completion of $F$ with respect to $v$, respectively, we have natural
isomorphisms
\[
K_q(F)/U^1K_q(F)\liso K_q(F_v)/U^1K_q(F_v) \liso K_q(\hat F_v)/U^1K_q(\hat F_v)
\]
for $q=1,2$. The usual approximation lemma for a finite set of discrete
valuations on a field shows the

\begin{lemma}\label{Kqapprox}
If\ $v_1,\ldots,v_s$ are finitely many different discrete valuations on a field $F$, then the diagonal
maps
\[
F^\times \to \bigoplus_{i=1}^s F_{v_i}^\times/U^1(F_{v_i}^\times)\quad \hbox{ and } \quad
K_2(F) \to \bigoplus_{i=1}^s K_2(F_{v_i})/U^1K_2(F_{v_i})
\]
are surjective.
\end{lemma}
Now let $Y$ be any integral scheme of finite type
over $\Spec(\Z)$. We call $Y$ horizontal if the generic point of $\Spec(\Z)$ is in the image
of the structural morphism. Otherwise we call $Y$ vertical, in which case it is a scheme of
finite type over a finite field. Let $Z \subset Y$ be a closed subscheme of positive
codimension. Consider the group
\[
\SK_1(Y,Z):=\coker \big(\bigoplus_{y \in (Y-Z)_1} K_{2}(k(y)) \to \bigoplus_{y \in Y_0}
K_1(k(y) \,\big).
\]
There is a natural surjection $\SK_1(Y,Z) \twoheadrightarrow \SK_1(Y)$ which is an
isomorphism if $Z$ is empty or has dimension zero. Any finite morphism $f:Y_1
\to Y_2$ with $f(Y_1-Z_1)\subset Y_2-Z_2$ induces a homomorphism $f_*: \SK_1(Y_1,Z_1) \to
\SK_1(Y_2,Z_2)$.

\medskip
If $Y$ is vertical, normal and proper we put $ \F=\Gamma(Y, \O_{Y})$. Note that $\F$ is the
algebraic closure of the prime field in $k(Y)$ and that $Y$ is a geometrically irreducible
variety over $\F$. The surjective norm maps $N: k(y)^\times \to \F^\times$ for $y \in Y_0$
give rise to a surjective norm map $\bigoplus_{y \in Y_0}k(y)^\times \to \F^\times$. It can
be easily deduced from the special case when $Y$ is a regular proper curve that the last map
annihilates the image of the boundary map $\bigoplus_{y \in Y_1} K_2(k(y)) \to \bigoplus_{y
\in Y_0}k(y)^\times$. Therefore we obtain a surjective norm map
\[
N: \SK_1(Y,Z) \to \F^\times.
\]

\pagebreak
\begin{prop} \label{sk1descr}
Let $Y$ be an integral scheme of finite type over $\Spec(\Z)$ and let $Z$ be a proper closed
subscheme. Then $\SK_1(Y,Z)$ is finite. More precisely:
\begin{itemize}
\item[\rm{(i)}] If\/ $Y$ is horizontal or vertical but not proper, then $\SK_1(Y,Z)=0$.
\item[\rm (ii)]
Assume that $Y$ is vertical and proper. Let $\tilde{Y}$ be the normalization of\/ $Y$ and let $Z_{\tilde{Y}}$ be preimage of $Z$ in $\tilde{Y}$. Then the natural map $\SK_1(\tilde{Y},Z_{\tilde{Y}})\to \SK_1(Y,Z)$ is surjective.
\item[\rm (iii)] If\/ $Y$ is vertical, proper and normal and\/ $\F=\Gamma(Y,\O_Y)$,
then the norm map
\[
N: \SK_1(Y,Z) \liso  \F^\times
\]
is an isomorphism.
\item[\rm (iv)] For any $y\in Y_0$, the natural homomorphism $k(y)^\times \to \SK_1(Y,Z)$ is
surjective.
\end{itemize}
\end{prop}

\begin{demo}{Proof:}
Two special cases of the proposition are classical, namely (i) if $Y$ is the spectrum of a
ring of integers in a number field and $Z=\varnothing$ (see \cite{BMS}) and (iii) for $Y$ a
smooth proper curve over a finite field and $Z=\varnothing$ (see \cite{Mo}).  We will deduce
the general assertion from this two special cases. First we note that the normalization
$\tilde{Y}\to Y$ is finite. Therefore the norm maps induce a well-defined homomorphism
$\SK_1(\tilde{Y},Z_{\tilde{Y}}) \to \SK_1(Y,Z)$, which is surjective, since $k(y)$ is a finite field for any closed point $y \in Y$. Now we proceed in several steps.

\noindent
1. Let $Y$ be horizontal and irreducible. In order to show the vanishing of $\SK_1(Y,Z)$, it
suffices to show that for every closed point $y \in Y$ the image of $k(y)^\times$ in
$\SK_1(Y,Z)$ is trivial. Let $y$ be an arbitrary closed point. Choose a closed irreducible
horizontal subscheme $Y'$ of $Y$ containing $y$ but not contained in $Z$. Putting $Z'=Y'\cap
Z$, the image of $k(y)^\times$ in $\SK_1(Y,Z)$ is contained in the image of the natural map
$\SK_1(Y',Z') \to \SK_1(Y,Z)$. By induction on the dimension, we may assume that $Y$ is of
dimension one. In this case we have $\SK_1(Y,Z)=\SK_1(Y)$ and $\tilde Y$ is the spectrum of the ring of $S$-integers in a number field. The surjection $\SK_1(\tilde {Y}) \to \SK_1(Y)$ together with the vanishing of $\SK_1(\tilde {Y})$ (see \cite{BMS}) shows that $\SK_1(Y)=0$.

\noindent
2. If $Y$ is affine and vertical, the same argument as in step 1, but now using Moore's
theorem \cite{Mo}, shows the vanishing of $\SK_1(Y,Z)$.

\noindent
3. Now assume that $Y$ is vertical, proper and normal and let $\F=\Gamma(Y,\O_Y)$. Then $Y$
is geometrically integral over $\Spec(\F)$. We first assume that $Y$ is projective. Let $H$
be a geometrically integral hypersurface section not contained in $Z$ (see \cite{Po} for its
existence) and $Z_H=Z\cap H$. We obtain an exact sequence
\[
\SK_1(\tilde{H},Z_{\tilde{H}}) \to \SK_1(Y,Z) \to \SK_1(Y-H,Z - Z_H),
\]
where $\tilde{H}$ is the normalization of $H$ and $Z_{\tilde{H}}$ is the preimage of $Z_H$ in $\tilde{H}$. Since $Y-H$ is affine, using step 2, the statement that the norm induces an isomorphism $N: \SK_1(Y,Z) \liso \F^\times$ can be reduced to the same statement for $(\tilde{H},Z_{{\tilde H}})$. By induction on the dimension, we may suppose that $Y$ is a curve. In this case the statement is just Moore's theorem \cite{Mo}. This solves the projective case. If $Y$ is proper but not necessarily projective over $\F$, there exists a birational morphism $f:Y' \to Y$ with $Y'$ normal and projective. Choose $W \subset Y$ containing $Z$ such that $f: Y'-f^{-1}(W) \liso Y-W$. Then we obtain a surjective homomorphism $\SK_1(Y',f^{-1}(W)) \to \SK_1(Y,Z)$, which shows the statement for $(Y,Z)$.

\noindent
4. Finally assume that $Y$ is vertical but not proper. In the same manner as in step~3 we may
suppose that $Y$ is normal and quasi-projective. Let $\bar Y$ be a projective
compactification. Choosing an irreducible hypersurface section $H$ not containing $Z$ but
containing at least one point of $\bar Y -Y$ and using step~2, we can proceed by induction
on the dimension. Finally, a non-proper curve is affine, showing $\SK_1(Y,Z)=0$.

\medskip\noindent
It remains to show (iv) if $Y$ is vertical and proper. If $Y$ is normal, the statement
follows from (iii). In the general situation, let $\tilde{y}$ be a point on the normalization $\tilde{Y}$ of $Y$ projecting to $y$. The statement follows from the commutative diagram
\diagram{cccl}
 k(\tilde{y})&\surjr{}& \SK_1(\tilde{Y},Z_{\tilde{Y}})\\
 \mapd{}&&\surjd{}\\
 k(y)& \to& \SK_1(Y,Z)&.
\enddiagram
\end{demo}

Now assume that $X$ is a connected noetherian excellent regular  scheme, pure of dimension
$d$, $D$ a divisor on $X$ and $Y=\supp(D)$.
We arrange the $E_1$-terms $E_1(X,D)$, $E_1(X|Y)$ and $E_1(Y^+)$ of the spectral sequences
introduced in the last section as a (up to sign convention) double complex in the following
way:

\diagram{ccccccc}
 0 & \to &E_1^{d,-d}(X,D)& \stackrel{\sim}{\to} &E_1^{d,-d}(X|Y)& \to &0\\
 \mapu{}&&\mapu{}&&\mapu{}&&\mapu{}\\
 E_1^{d-1,-d-1}(Y^+) & \to& E_1^{d-1,-d}(X,D)& \to &E_1^{d-1,-d}(X|Y)& \stackrel{f}{\to} &E_1^{d-1,-d}(Y^+)\\
 \mapu{}&&\mapu{}&&\mapu{g}&&\mapu{}\\
 E_1^{d-2,-d-1}(Y^+) & \to& E_1^{d-2,-d}(X,D)& \to &E_1^{d-2,-d}(X|Y)&
 \stackrel{h}{\to} &E_1^{d-2,-d}(Y^+)
 \enddiagram
The rows are exact. The vertical maps are the differentials of the various spectral
sequences. In the notation of proposition~\ref{chow0descrklein}, $R_{X,D}=\ker(f)$. The zeros
in the diagram are due to dimension reasons. A diagram chase shows the

\begin{lemma} \label{fincrit}
There exists a natural exact sequence
\[
\im(f)/\im(f\circ g) \lang \CH_0(X,D) \lang E_2^{d,-d}(X|Y) \lang 0.
\]
\end{lemma}
In order to analyze the exact sequence of lemma~\ref{fincrit}, we first deduce some
approximation results.  Note that the maps of the next lemma are not the maps from the above diagram, unless $Y$ is irreducible.

\begin{lemma} \label{surjlemma}
The natural maps
\begin{description}
\item{\rm (i)} $E_1^{d-2,-d}(X|Y) \lang E_1^{d-2,-d}(Y)$
\item{\rm (ii)} $E_1^{d-1,-d}(X|Y) \lang E_1^{d-1,-d}(Y)$
\item{\rm (iii)} $E_1^{d-2,-d-1}(X|Y) \lang E_1^{d-2,-d-1}(Y)$
\end{description}
 are surjective.
\end{lemma}

\begin{demo}{Proof:}
Let $y_1 \in Y^{d-2}$ and choose a point $x \in (X-Y)^{d-2}$ such that $y_1$ is a regular
point of the closure $W$ of ${\{x\}}$ in $X$.  Let $F=k(W)$ be the function field and $W'$
the normalization of $W$.  Let $y'_1, \ldots y'_s \in (W')^1$ be the finitely many points of
$(W')^1$ lying over point in $Y^{d-2}$. We choose the ordering such that $y'_1$ is the unique
point over $y_1$ (by assumption, $y_1$ is a regular point on $W$). Each $y' \in
\{y'_1,\ldots,y'_s\}$ defines a discrete valuation on $F$. Let $M$ be the kernel of the
(multi) tame symbol $K_2(F) \to \bigoplus_{i=1}^s k(y'_i)^\times$. By
proposition~\ref{ausrech}, each $m\in M$ defines an element in
$E_1^{d-2,-d}(X|Y)$. We obtain a commutative diagram with exact lines

\diagram{ccccccccl}
 0&\to&M&\to&K_2(F)&\to&\ds \bigoplus_{i=1}^s k(y'_i)^\times&\to&0 \\
&& \mapd{}&&\mapd{}&&\eqd\\
0&\to&\ds\bigoplus_{i=1}^s K_2(k(y'_i))&\to&\ds\bigoplus_{i=1}^s
K_2(F_{y'_i})/U^1K_2(F_{y'_i}) &\to & \ds\bigoplus_{i=1}^s k(y'_i)^\times &\to& 0\,.
\enddiagram
By lemma~\ref{Kqapprox}, the middle vertical arrow is surjective and therefore we find an $m
\in M$ mapping to an arbitrarily chosen $\alpha \in K_2(k(y_1))=K_2(k(\tilde{y}_1))$ and to
zero in all other components. This shows (i). The proofs of (ii) and (iii) are similar and
left to the reader.
\end{demo}

\begin{prop} \label{e2iso}
The natural map
\[
E_2^{d,-d}(X|Y) \lang E_2^{d,-d}(X) = \CH^d (X)
\]
is an isomorphism. If $X$ is flat and of finite type over $\Spec(\Z)$, then both groups are
finite.
\end{prop}

\begin{demo}{Proof:}
By proposition~\ref{ausrech}, $K_1(\M^{d-1/d}_{X-Y})\mapr{\delta} K_0(\M_Y^{d-1/d})$ is
surjective. Therefore, any class in $\CH^d(X)$ can be represented by a cycle with support in
$X-Y$, which shows the surjectivity. Next consider the commutative diagram
\diagram{ccccc}
E_1^{d,-d}(X|Y)&\stackrel{a}{\hookrightarrow}&E_1^{d,-d}(X)&\twoheadrightarrow
& E_1^{d-1,-d-1}(Y)\\
\mapu{}&&\mapu{b}&&\mapu{}\\
E_1^{d-1,-d}(X|Y)&\hookrightarrow &E_1^{d-1,-d}(X)&\twoheadrightarrow & E_1^{d-2,-d-1}(Y)\\
\mapu{}&&\mapu{}\\
E_1^{d-2,-d}(X|Y)&\to&E_1^{d-2,-d}(X).
\enddiagram
The lines and rows are complexes and, by proposition~\ref{ausrech}, the upper line is exact.
If $\alpha\in E_1^{d,-d}(X|Y)$ represents an element in $\ker(E_2^{d,-d}(X|Y)\to
E_2^{d,-d}(X))  $, we find a $\beta \in E_1^{d-1,-d}(X)$ with $a(\alpha)=b(\beta)$. By
lemma~\ref{surjlemma} (iii), the composite map $E_1^{d-2,-d}(X|Y)\to E_1^{d-2,-d-1}(Y)$ is
surjective. Therefore, we may choose $\beta \in E_1^{d,-d}(X-Y)\subset E_1^{d-1,d}(X)$. By
Proposition~\ref{ausrech}, we obtain $\beta \in E_1^{d-1,-d}(X|Y)$, which shows the
injectivity. Finally, if $X$ is flat and of finite type over $\Spec(\Z)$ then $\CH^d(X)$ is
finite by \cite{K-S3}, theorem~6.1.
\end{demo}

Now we assume that $X$ is of finite type over $\Spec(\Z)$ and let, for $i=1,\ldots, r$,   $Z_i= \bigcup_{j\neq i} D_j\cap D_i \subset
D_i$. Put $Z=\tcup_{i=1}^r Z_i \subset Y^+$. Let $f$, $g$ and $h$ be the maps in the diagram
before lemma~\ref{fincrit}.

\begin{lemma}\label{approx3}
\begin{description}
\item{\rm (i)} $\im(f) \supset \bigoplus_{y \in (Y^+-Z)_0} k(y)^\times$.
\item{\rm (ii)} $\im (h) \supset \bigoplus_{y \in (Y^+-Z)_1} K_2(k(y))$.
\end{description}
\end{lemma}

\begin{demo}{Proof:}
This follows from lemma~\ref{surjlemma} (i) and (ii).
\end{demo}

Let $\tilde{Y}$ be the normalization of\/ $Y$.
\begin{prop} We have a natural surjection
\[
\SK_1(\tilde{Y}) \surjr{} \im(f)/\im(f\circ g).
\]
\end{prop}

\begin{demo}{Proof:}
First note that $\tilde{Y}$ is also the normalization of $Y^+$. Let  $Z_{\tilde{Y}}$ be the preimage of\/ $Z$ in
$\tilde{Y}$. By lemma~\ref{approx3} (ii), we obtain surjections
\[
\SK_1(Y^+,Z) \surjr{} E_1^{d-1,-d}(Y^+)/\im (f\circ g) \surjr{} \SK_1(Y^+).
\]
Hence, by proposition~\ref{sk1descr} (iv), and lemma~\ref{approx3} (i), the natural inclusion
map
\[
\im (f)/\im (f\circ g) \lang E_1^{d-1,-d}(Y^+)/\im (f\circ g)
\]
is an isomorphism. By proposition~\ref{sk1descr} (ii), the natural map
\[
\SK_1(\tilde{Y},Z_{\tilde{Y}})\lang \SK_1(Y^+,Z)
\]
is surjective. Finally, by proposition~\ref{sk1descr} (iii), $\SK_1(\tilde{Y},Z_{\tilde{Y}})\cong \SK_1(\tilde{Y})$.
\end{demo}

 Using propositions~\ref{e2iso} and  \ref{sk1descr}, lemma~\ref{fincrit} implies the

\begin{theorem} \label{chrelfinite}
Let $X$ be a regular connected scheme of finite type over $\Spec(\Z)$ and let $D=D_1+\cdots
+D_r$ be a sum of prime divisors on $X$. Denoting the normalization of $D_i$ by
$\tilde{D}_i$ for $i=1,\ldots, r$,  we have a natural exact sequence
\[
\bigoplus_{i=1}^r \SK_1(\tilde{D}_i) \lang \CH_0(X,D) \lang \CH_0(X) \lang 0.
\]
In particular, if $X$ is flat over $\Spec(\Z)$, the group $\CH_0(X,D)$ is finite.
\end{theorem}

Finally, we enlarge the divisor. Let $D_1, \ldots , D_{r+s}$ be prime divisors on $X$, $D=
D_1+\cdots+D_r$ and $D'=D_{1}+\cdots+D_{r+s}$. Let $Y$ and $Y'$ be the support of $D$ and
$D'$, respectively, and let $\CH_0^{D'}(X,D)$ be the quotient of $\ds\bigoplus_{x\in
(X-Y')_0}\!\!\Z\;$ by the image of the group
\[
R_{X,D,D'}\stackrel{\rm df}{=}\ker \Big( \ker \big(\bigoplus_{x \in X_1-Y'} k(x)^\times \to
\bigoplus_{y \in Y'_0} \Z \big) \lang  \bigoplus_{z \in {(Y^+)}_0} k(z)^\times \Big)
\]
under the divisor map.

\begin{prop} \label{horver1}
Let $X$ be a regular scheme, flat and of finite type over $\Spec(\Z)$ and let\break $D_1, \ldots ,
D_{r+s}$ be pairwise different prime divisors on $X$. Let $D= D_1+\cdots+D_r$,
$D'=D_{1}+\cdots+D_{r+s}$ and denote the normalization of $D_i$ by $\tilde{D}_i$,
$i=1,\ldots, r+s$. Then the following holds.
\begin{itemize}
\item[\rm (i)] $\CH_0^{D'}(X,D) \liso \CH_0(X,D)$.
\item[\rm (ii)]
We have a natural exact sequence
\[
\bigoplus_{i=r+1}^{r+s} \SK_1(\tilde{D}_i) \lang \CH_0(X,D') \lang \CH_0(X,D) \lang 0.
\]
In particular, if $D_{r+1},\ldots,D_{r+s}$ are horizontal, then $\CH_0(X,D') \to \CH_0(X,D)$
is an isomorphism.
\end{itemize}
\end{prop}

\begin{demo}{Proof:} By the definition of the occurring objects, we have natural maps
\[
\CH_0(X,D') \mapr{p} \CH_0^{D'}(X,D) \mapr{q} \CH(X,D)
\]
The map $\bigoplus_{i=r+1}^{r+s} \SK_1(\tilde{D}_i) \lang \CH_0(X,D')$ is the composite map
\[
\bigoplus_{i=r+1}^{r+s} \SK_1(\tilde{D}_i) \mapr{\incl} \bigoplus_{i=1}^{r+s}
\SK_1(\tilde{D}_i) \lang \CH_0(X,D')
\]
of the natural inclusion with the map of theorem~\ref{chrelfinite}. The commutative diagram
\diagram{ccccccc}
\bigoplus_{i=r+1}^{r+s} \SK_1(\tilde{D}_i)\\
\injd{}\\
\bigoplus_{i=1}^{r+s} \SK_1(\tilde{D}_i)&\to&\CH_0(X,D')&\to&\CH_0(X)&\to &0\\
\surjd{}&&\mapd{q\circ p}&&\eqd\\
\bigoplus_{i=1}^{r} \SK_1(\tilde{D}_i)&\to&\CH_0(X,D)&\to&\CH_0(X)&\to &0
\enddiagram
shows that the sequence in (ii) is a complex and that $q\circ p$ is surjective.  We consider
two copies of the double complex before lemma~\ref{fincrit}: one for $(X,D)$
and one for $(X,D')$, and the natural map between them. Assume that the class of a
zero-cycle $\sum \alpha_i P_i$, $P_i \in X-Y'$, is zero in $\CH_0(X,D)$. By
proposition~\ref{chow0descrklein}, $\sum \alpha_i P_i$, is the image of an element $a\in
R_{X,D}$ under the divisor map. Changing $a$ by the image of an appropriate element in
$E_1^{d-2,-d}(X|Y)$ (use approximation as in the proof of lemma~\ref{surjlemma}), we may
suppose that the components of $a$ at all points of $Y'-Y$ are trivial, i.e.
\[
a\in \ker\Big(\ker \Big( \bigoplus_{x \in (X-Y')_1} k(x)^\times \lang \bigoplus_{y\in Y_0} \Z
\Big)\lang \bigoplus_{z \in Y^+}k(z)^\times\Big).
\]
Since all (non-zero) components of $\div(a)=\sum n_i P_i$ are in $X-Y'$, we conclude $a\in
R_{X,D,D'}\subset E_1^{d-1,-d}(X|Y')$, proving (i).
With respect to the map $f': E_1^{d-1,-d}(X|Y'){\to} E_1^{d-1,-d}((Y')^+)$, we obtain
$$
f'(a) \in \bigoplus_{z \in ((Y')^+-Y^+)_0}k(z)^\times \subset E_1^{d-1,-d}((Y')^+).
$$
Therefore the class in $\ds\bigoplus_{i=r+1}^{r+s} \SK_1(\tilde{D}_i)$ of any preimage of
$f'(a)$ in $\ds\bigoplus_{i=r+1}^{r+s} \bigoplus_{z \in (\tilde D_i)_0}k(z)^\times$ maps to
the class of $\sum \alpha_i P_i$ in $\CH_0(X,D')$. This completes the proof.
\end{demo}

\section{Tame coverings} \label{tamesect}

Coverings of a regular scheme which are tamely ramified along a normal crossing divisor have
been studied in \cite{SGA1}, \cite{G-M}. A naive extension of the (valuation theoretic)
definition of tame ramification in the normal crossing case proves to be not useful in the
general situation. For example, it would not be stable under base change (cf. \cite{S1},
Example 1.3). A definition of tameness in the general situation was given in \cite{S1} and
it was shown there that it coincides with the previous one in the normal crossing case. Here
we restrict our attention to  tamely ramified extensions of normal schemes, where one can
define tameness using inertia groups.

\medskip
 Let $X$ be a connected, noetherian, normal scheme and let $K$ be its
function field. Let $L|K$ be a finite Galois extension and let $X_L$ be the normalization of
$X$ in $L$. Furthermore, let $U\subset X$ be an open subscheme such that $U_L \to U$ is
\'{e}tale. Put $Y=X-U$.

\begin{defi} \label{tamedefprop}
We say that $U_L \to U$ is tamely ramified along $Y$ if  the inertia subgroup (cf.\ {\rm
\cite{Bo}, V,\S2.2}) of every point $y \in Y$ in $G(L|K)$ is of order prime to the
characteristic of\/ $k(y)$.
\end{defi}

If $X$ is regular, $D=D_1 + \cdots +D_r$ a divisor with normal crossings on $X$ and\/
$Y=\supp(D)$, then $U_L \to U$ is tamely ramified along $Y$ if and only if $L|K$ is tamely
ramified at the discrete valuations $v_1,\ldots,v_r$ of $K$ corresponding to $D_1,\ldots,
D_r$ (see \cite{S1}, prop.1.14).

\medskip
We say that $X_L \to X $ is tame along a divisor $D$, if, putting $Y=\supp(D)$ and $U=X-Y$,
the induced morphism $U_L\to U$ is \'{e}tale and tame along $Y$. For any base point $* \in
U=X-Y$, the tame fundamental group $\pi_1^t(X,D,*)$ is the unique quotient of $\pi_1(U,*)$
which classifies finite connected \'{e}tale coverings $\tilde{U}$ of $U$ with at most tame
ramification along $D$.

\medskip
If $X$ is of finite type over $\Spec(\Z)$, we denote  by $\tilde \pi_1^t(X,D,*)^{ab}$ the
unique quotient of $\pi_1^t(X,D,*)$ which classifies abelian finite \'{e}tale coverings of
$U=X-Y$ which are at most tamely ramified along $Y=\supp(D)$ and in which every real (i.e.\
$\R$-valued) point of $U$ splits completely. Fundamental groups with respect to different
base points are isomorphic, the isomorphism being canonical up to inner automorphisms. As
usual, we omit the base point as long as we deal with abelian quotients.

\begin{prop} \label{hotver2}
Let $X$ be a connected regular scheme, proper and flat over $\Spec(\Z)$ and let $D$ be a
divisor on $X$. Then for any {\em horizontal} prime divisor  $D'$ on $X$ the natural map
\[
\tilde\pi_1^t(X,D+D')^\ab \lang \tilde \pi_1^t(X,D)^\ab
\]
is an isomorphism.
\end{prop}

\begin{demo}{Proof:}
The map in question is obviously surjective. To show injectivity, we first show that every
connected cyclic \'{e}tale covering $\tilde V \to V:=X-(Y\cup D')$ of $p$-power degree ($p$ an
arbitrary prime number) which is tamely ramified along $Y \cup D'$ extends to an \'{e}tale
covering of $U:=X-Y$. Let $K$ and $L$ be the function fields of $V$ and of $\tilde V$,
respectively. Since $D'$ is horizontal, the map $D' \to \Spec(\Z)$ is proper and dominant,
hence surjective. Thus $D'$ contains a point of residue characteristic $p$, which has a
trivial inertia group in $L|K$ since $\tilde V \to V$ is tame along $D'$. Hence $L|K$ is
unramified at the discrete valuation associated to $D'$. The theorem on the purity of the
branch locus then shows that the normalization $\tilde U$ of $U$ in $L$ is \'{e}tale over $U$.
This shows that $\pi_1^t(X,D+D')^\ab \to \pi_1^t(X,D)^\ab$ is an isomorphism. Finally, we
have to deal with the real places.  By \cite{Sa}, lemma 4.9 (iii), the subset of points in
$U(\R)$ which split completely in $\tilde U$\/ is (norm) closed and open in $U(\R)$. $U(\R)$
is a real manifold and $V(\R)$ is (norm) dense in $U(\R)$. Hence, if all points of $V(\R)$
split completely, in $\tilde V$, then all points of $U(\R)$ split completely in $\tilde U$.
Therefore the result extends to the modified fundamental groups.
\end{demo}

The following finiteness result is
a special case of (\cite{S1}, Theorem 2).

\begin{theorem}\label{endlichkeit}
Let\/ $X$ be a connected, regular, proper and flat scheme over $\Spec(\Z)$ and
let\/ $D$ be a divisor on\/ $X$. Then the group $\tilde \pi_1^t(X,D)^{ab}$ is finite.
\end{theorem}

We conclude this section with a lemma which will be needed later on.

\begin{lemma} \label{tamebyrank1}
Let $X$ be a connected regular scheme, flat and of finite type over $\Spec(\Z)$. Let $D=D_1 +
\cdots +D_r$ be a sum of\/ {\em vertical} prime divisors on $X$ (not necessarily with normal
crossings) and put $U=X-\supp(D)$. Let $v_1,\ldots,v_r$  be the discrete valuations of the
function field $k(U)=k(X)$ which are associated with $D_1,\ldots,D_r$. Then a finite abelian
\'{e}tale covering $\tilde U \to U$ is tamely ramified along $D$ if and only if the extension
$k(\tilde U)|k(U)$ is tamely ramified at $v_1,\ldots,v_r$.
\end{lemma}

\begin{demo}{Proof:}
In order to show the nontrivial implication, assume that $k(\tilde U)|k(U)$ is tamely
ramified at $v_1,\ldots,v_r$. We may assume that $\tilde U \to U$ is cyclic of prime power
order, say of order $p^n$. After reordering, we may assume that $D_1,\ldots,D_s$ lie over
$p$ and that $D_{s+1},\ldots,D_r$ lie over prime numbers different to $p$. Put
$V=X-\supp(D_{s+1}+\cdots +D_r)$.  Since $v_1,\ldots, v_s$ are tamely ramified, hence
unramified in $k(\tilde U)$, the theorem on the purity of the branch locus shows that the
normalization $\tilde V$ of $V$ in $k(\tilde U)$ lies \'{e}tale over $V$. Since every point of
$\supp(D)$ with residue characteristic $p$  lies on $V$, this shows that $\tilde U \to U$ is
tamely ramified along $D$.
\end{demo}

\section{Some local computations} \label{recexsect}
Let $X$ be an excellent regular scheme, $i:Z \hookrightarrow X$ a closed regular subscheme
of pure codimension~$c$, $p$ a prime number with $\frac{1}{p}\in \O_X$ and $n \geq 1$ an integer. All  cohomology groups occurring in this section are taken with respect to the {\'e}tale topology. The cup-product with the fundamental class (cf.\ \cite{De}) $\cl(Z)\in
H^{2c}_D(X,\Z/p^n\Z(c))$ induces the so-called Gysin map
\[
\Gys: i^*(\Z/p^n\Z) \lang Ri^!(\Z/p^n\Z)(c)[2c].
\]
It is a special case of ``purity for \'{e}tale cohomology'' (proved by Gabber, cf.\ \cite{Fu})
that $\Gys$ is an isomorphism. The following proposition follows from this result in a straightforward manner.

\begin{prop}\label{lokalekoho}
Let $X$ be an excellent regular scheme, $i:Y \hookrightarrow X$ a closed subscheme
of codimension $\geq c$, $p$ a prime number with $\frac{1}{p}\in \O_X$ and $n \geq 1$ an integer. Furthermore, let\/ $Y^\sing$ be the singular locus of\/ $Y$ and $Y^\ns=Y-Y^\sing$. Then, for any $i\in \Z$, we have
\[
H^r_Y(X, \Z /p^n\Z (i))= 0 \quad \hbox{ for } r < 2c.
\]
If\/ $Y$ is of pure codimension $c$ in $X$, then we have natural isomorphisms
\[
H^{2c}_Y(X, \Z /p^n\Z (i))\liso H^{2c}_{Y^\ns}(X-Y^\sing, \Z /p^n\Z (i))\rightiso H^0(Y^\ns, \Z /p^n\Z (i-c)).
\]
The first map is the natural restriction homomorphism and the second one is the cup-product with the fundamental class $c_{Y^\ns} \in  H^{2c}_{Y^\ns}(X-Y^\sing, \Z /p^n\Z (c))$.
\end{prop}

\begin{demo}{Proof:} If $Y$ is regular, this follows from purity. In order to deal with the general case, we proceed by (decreasing) induction on the codimension and use the long exact sequence
\[
\cdots \to H^r_{Y^\sing}(X,\Z /p^n\Z (i)) \to H^r_{Y} (X, \Z /p^n\Z (i)) \to H^r_{Y-Y^\sing}(X-Y^\sing, \Z /p^n\Z (i)) \to \cdots\; .
\]
\end{demo}

\bigskip
Now let $D=D_1+\cdots+D_r$ be the sum of distinct vertical irreducible divisors on $X$ (i.e.\ the $D_i$'s are schemes over fields), and we put $Y=\supp(D)$. Let $p$ be a prime number and let
$Y_p$ be the inverse image of $\{p \}\subset \Spec(\Z)$ under the natural projection $Y\to
\Spec(\Z)$. An abelian \'{e}tale covering of $X-Y$ of $p$-power degree which is tame along $Y$ extends to an \'{e}tale covering of $X-(Y-Y_p)$. Since all $D_i$ are vertical, tameness of such a covering along $Y-Y_p$ is a void condition. Therefore we have a natural isomorphism
$\pi_1^t(X,D)^\ab(p) \liso \pi_1(X-(Y-Y_p))^\ab(p)$. Using excision and proposition~\ref{lokalekoho}, we have
\[
H^2_{Y-Y_p}(X,\Q_p/\Z_p)=H^2_{Y-Y_p}(X-Y_p,\Q_p/\Z_p)=H^0(Y-(Y_p \cup Y^\sing), \Q_p/\Z_p(-1))
\]
and therefore an exact sequence
\[
H^0(X-(Y_p \cup Y^\sing), \Q_p/\Z_p(-1))^\vee \to \pi_1^t(X,D)^\ab(p)\to \pi_1(X)^\ab(p) \to 0.
\]
For a field $K$, we use the notation
\[
H^0(K,\Q/\Z(-1))= \dirlim_{n} \mu_n(K)^\vee,
\]
where $\scriptstyle \vee$ denotes Pontryagin dual. In particular, the char$(K)$-part of the above group is trivial. We obtain the

\begin{lemma}\label{Gysin}
Let $D=D_1+\cdots+D_r$ be the sum of distinct vertical irreducible divisors on an excellent regular
scheme $X$.  Then we obtain a natural exact sequence
\[
\bigoplus_{i=1}^r H^0(k(D_i),\Q/\Z(-1))^\vee \lang \pi_1^t(X,D)^\ab \lang \pi_1(X)^\ab \lang
0.
\]
\end{lemma}

Next we introduce a technical condition.

\begin{defi}\label{condmu}
Let $Y$ be an integral scheme of finite type over a finite field\/ $\F$. We say that $Y$ satisfies {\boldmath\bf condition $(\mu)$} if for every natural number $n$, prime to the characteristic of\/ $\F$, the natural map 
\[
H^0(Y,\Z/n\Z(-1)) \injr{} H^0(k(Y), \Z/n\Z(-1))
\]
is an isomorphism.
\end{defi} 
Condition $(\mu)$ is satisfied if $Y$ is normal.  If $Y$ is proper,  condition $(\mu)$ is equivalent to the condition that $\Gamma(Y,\O_Y)$ coincides with the algebraic closure of $\F$ in $k(Y)$.

\medskip\noindent
{\it Example.} If $p$ is a prime number, $p\equiv 3 \bmod 4$, then $Y= \Spec(\F_p[x,y]/x^2+y^2)$ does not satisfy condition $(\mu)$. Indeed, $k(Y)$ is the rational function field in one variable over $\F_{p^2}$.

\medskip
Now let $\Z_p^h$ be the henselization of $\Z$ at a prime $p$. Let $X\to \Spec(\Z_p^h)$ be a regular, proper scheme of finite type and $X_\eta$ and $X_s$ its generic and special fibre. Let $D=D_1+\cdots+D_r$ be the sum of distinct vertical irreducible divisors on $X$ and $Y=\supp(D)$. Local class field theory induces a natural map 
\[
\tau: \bigoplus_{x\in {(X_\eta)}_0} k(x)^\times \lang \bigoplus_{x\in {(X_\eta)}_0} \pi_1(x)^\ab \lang \pi_1^t(X,D)^\ab
\]
and our next aim is to detect a natural subgroup of its kernel. This will be the key point in the proof of the existence of the global reciprocity map in the next section. Its exactly here, where we have to impose the above condition $(\mu)$ on the divisors $D_i$. 
We consider the group
\[
U^t_D(X)\stackrel{\mathrm{df}}{=}\ker \Big( \ker \big(
\bigoplus_{(X_\eta)_0}
k(x)^\times \to  \bigoplus_{y\in (X_s)_0} \Z\big) \mapr{\phi}
  \bigoplus_{y\in (Y^+)_0} k(y)^\times \Big),
\]
where $\phi$ is the restriction of the natural map $E_1^{d-1,-d}(X|Y) \to E_1^{d-1,d}(Y^+)$, $d=\dim X$, of the spectral sequence for relative $G$-theory from section~\ref{relk} to the subgroup of those elements which have nontrivial components only at points in the generic fibre. $U^t_D(X)$ is nothing else but the group $R_{X,D,D'}$ defined before proposition~\ref{horver1} with $D'$ the sum of all vertical divisors of $X$. A point $x\in (X_\eta)_0$ specializes to a unique point $y\in (X_s)_0$. Therefore, putting
\[
U_D^t(y)=  
 \ker \Big( \ker \big(
\bigoplus_{
\renewcommand{\arraystretch}{0.8}\begin{array}{c} \sst x \in (X_\eta)_0 \\
\sst x \to y
\end{array}
\renewcommand{\arraystretch}{1}}
k(x)^\times \to  \Z\big) \mapr{\phi}
  \bigoplus_{ D_i\ni y} k(y)^\times \Big),
\]
we obtain a decomposition
\[
U^t_D(X)= \bigoplus_{y\in (X_s)_0} U^t_D(y).
\]

\begin{prop}\label{utd}  Let $X\to \Spec(\Z_p^h)$ be a regular, proper scheme of finite type and let $D=D_1+\cdots+D_r$ be the sum of distinct vertical irreducible divisors on $X$. If $D_1,\ldots,D_r$ satisfy condition $(\mu)$, then the natural map
\[
\tau: \bigoplus_{x \in (X_\eta)_0} k(x)^\times  \lang \pi_1^t(X,D)^\ab
\]
annihilates the subgroup $U^t_D(X)$.
\end{prop} 

\begin{demo}{Proof:} We consider the $l$-part for all prime numbers $l$ separately. The easiest case is $l=p$. Indeed, in the  commutative diagram (cf.\ A.3)
\diagram{ccccccc}
\ds \bigoplus_{x \in {(X_\eta)}_0} k(x)^\times &\lang&\SK_1(X_\eta)&\mapr{\tau}&\pi_1(X_\eta)^\ab &\lang& \pi_1^t(X,D)^\ab\\
\mapd{}&&\mapd{\partial}&&\mapd{\sp}&&\mapd{}\\
\ds \bigoplus_{y\in (X_s)_0} \Z &\lang&\SK_0(X_s)&\mapr{\phi}&\pi_1(X_s)^\ab&\rightiso& \pi_1(X)^\ab\;,
\enddiagram
 the right vertical arrow is an isomorphism on $p$ parts (by lemma~\ref{Gysin}). 
 
 \medskip\noindent
For the prime-to-$p$ part, we first note that it suffices to consider $U^t_D(y)$ for all $y \in (X_s)_0$. We fix once and for all such a $y$ and renumber the $D_i$ such that $D_i \ni y$ for $i=1,\ldots s$ and $D_i \not\ni y$ for $i=s+1,\ldots, r$. Now let $n$ be a natural number with $(n,p)=1$. Let $x_1,\ldots, x_m \in (X_\eta)_0$ be points that specialize to $y$ and let $W$ be the closure of $\{x_1,\ldots,x_m\}$ in $X$. $W$ is a one dimensional, local henselian scheme with finite residue field. In particular, $H^2(W,\Z/n\Z)=0$ and we obtain a short exact sequence
\[
0\lang H^1(W,\Z/n\Z) \lang \bigoplus_{j=1}^m H^1(k(x_i),\Z/n\Z) \mapr{\delta} H^2_y(W,\Z/n\Z)\lang 0. 
\]
Now consider the commutative diagram
\diagram{ccc}
\ds\bigoplus_{j=1}^r H^0(k(D_j),\Z/n\Z(-1))& \liso & H^2_{Y^\ns} (X-Y^\sing, \Z/n\Z)\\
\inju{\alpha}&&\mapu{\wr}\\
\ds\bigoplus_{j=1}^r H^0(D_i,\Z/n\Z(-1)) &\mapr{\sum (? \cup c_{D_j})} & H^2_{Y} (X, \Z/n\Z)\\
\mapd{\mathop{\oplus}\limits_{j=1}^s i^*_j}&&\mapd{i^*}\\
\bigoplus_{j=1}^s H^0(y,\Z/n\Z(-1)) &\mapr{\sum (? \cup i_j^*c_{D_j})} & H_y^2(W,\Z/n\Z)
\enddiagram
Here $i:(W,y)\hookrightarrow (X,Y)$ and $i_j: y \hookrightarrow D_j$, $j=1,\ldots,s$, are the inclusion morphisms and $c_{D_j} \in H^2_{D_j}(X,\Z/n\Z(1))$, $j=1,\ldots,r$, is the fundamental class of the divisor $D_j$ in $X$ (\cite{De}, 2.1). The horizontal cup-products are defined in \cite{De}, 1.2. Since the $D_j$ satisfy condition $(\mu)$, the inclusion map $\alpha$ is an isomorphism, and hence the same is true for the middle horizontal arrow $ \sum(? \cup c_{D_j})$.

Dualizing, and using the well-known duality theorems for Galois cohomology of finite and of one-dimensional henselian fields, we obtain a commutative diagram
\diagram{ccccc}
&&\hspace*{-2cm}\ker\big(\ds\bigoplus_{i=1}^m H^1(k(x_i),\Z/n\Z(1)) \to \Z/n\Z\big) & \lang &\ds\bigoplus_{i=1}^m H^1(k(x_i),\Z/n\Z(1))\\
&&\mapd{\wr}&&\mapd{\wr}\\
\ds\bigoplus_{j=1}^sH^1(k(y),\Z/n\Z(1))&\mapl{}&H^2_y(W,\Z/n\Z)^\vee&\lang&H^1(W-y,\Z/n\Z)^\vee\\
\mapd{}&&\mapd{}&&\mapd{}\\
\ds\bigoplus_{j=1}^sH^0(D_j,\Z/n\Z(-1))^\vee&\mapl{\sim}& H^2_{Y}(X,\Z/n\Z)^\vee&\lang & H^1(X-Y,\Z/n\Z)^\vee
\enddiagram

In order to show the proposition, it therefore remains to show  the commutativity of the following diagram
\diagram{ccc}
\ker \big( \ds\bigoplus_{i=1}^s k(x_i)^\times \lang \Z) &\lang &\ker\big(\ds\bigoplus_{i=1}^m H^1(k(x_i),\Z/n\Z(1)) \to \Z/n\Z\big)\\
\mapd{\phi}&&\mapd{\psi}\\
\ds\bigoplus_{j=1}^s k(y)^\times &\lang& \ds\bigoplus_{j=1}^s H^1(k(y),\Z/n\Z(1)),
\enddiagram
in which the horizontal maps are induced by the respective Kummer sequences, $\phi$ is the map from the definition of $U^t_D(y)$ and $\psi$ is the map induced by the last diagram. We gave a description of $\phi$ in section~\ref{explsect}. It remains to describe $\psi$.

Choose a  a local parameter $\pi_j$ defining $D_j$ in a neighborhood of $y$ in $X$. Then $\pi_j$ defines an element in $H^0(W-y,\G_m)$ and, denoting its image in $H^1(W-y,\Z/n\Z(1))$ by $\bar \pi_j$, the pull-back of the fundamental class $i_j^*(c_{D_j}) \in H^2_y(W,\Z/n(1))$ is the image of $\bar \pi_j$ under the natural boundary map $H^1(W-y,\Z/n\Z(1)) \to H^2_y(W,\Z/n\Z(1))$ (see \cite{De}, 2.1.1.1). We therefore obtain a commutative diagram 
\diagram{ccc}
H^0(y,\Z/n\Z (-1)) & \mapr{?\cup i_j^*c_{D_j}} & H^2_y(W,\Z/n\Z)\\
\mapu{\wr}&&\mapu{\delta}\\
H^0(W, \Z/n\Z(-1))& \mapr{?\cup \bar \pi_j} & H^1(W-y,\Z/n\Z)
\enddiagram

Next we recall the boundary map in {\'e}tale homology (see \cite{J-S}, \cite{B-O})
\[
b: H^i(W-y, \Z/n\Z (j)) \lang H^{i-1}(y,\Z/n\Z(j-1)), 
\]
which, denoting the structure map by $f: W \to \Spec(\Z_p^h)$, is defined by 
\[
H^i(W-y,\Z/n\Z(j))=H^i(W-y,f^!\Z/n\Z(j)) \mapr{\delta} H^{i+1}(y, i^!f^!\Z/n\Z(j))= H^{i-1}(y,\Z/n\Z(j-1)).
\]
We have the following commutative diagram relating $b$ with the symbol map $h: K_*(F) \to H^*(F,\Z/n\Z(*))$, $*=1,2$, from $K$-theory to Galois cohomology 
\[
\begin{array}{ccc}
\bigoplus_{i=1}^m K_2(k(x_i))& \mapr{h} & H^2(W-y, \Z/n\Z(2))\\
\mapd{\partial}& &\mapd{b}\\
k(y)^\times &\mapr{h}& H^1(k(y),\Z/n(1)),
\end{array}
\]
where $\partial$ is the boundary map of the Quillen spectral sequence for $W$ (see \cite{Ka}, Lemma 1.4 (1) for the case when $W$ is regular, the general case follows easily from that). Furthermore, we have the following commutative diagram, relating the Tate-duality for $W-y$ with that of $y$
\diagram{ccc}
H^2(W-y,\Z/n\Z(2)) &\cong& H^0(W-y, \Z/n\Z(-1))^\vee\\
\mapd{b}&&\mapd{\can}\\
H^1(y,\Z/n\Z(1))& \cong &H^0(y,\Z/n\Z(-1))^\vee\,.
\enddiagram 
Here $\can$ is the dual to the natural map 
\[
H^0(y,\Z/n\Z(-1))\cong H^0(W,\Z/n\Z(-1))\to H^0(W-y,\Z/n\Z(-1)).
\]
Summing up, the $j$-component of $\psi$ is the restriction of the composite map
\[
\bigoplus_{i=1}^m H^1(k(x_i),\Z/n\Z(1)) \mapr{?\cup \bar \pi_j} 
\bigoplus_{i=1}^m H^2(k(x_i),\Z/n\Z(2)) \mapr{b} H^1(k(y),\Z/n\Z(1))
\]
to the subgroup $\ker\big(\ds\bigoplus_{i=1}^m H^1(k(x_i),\Z/n\Z(1)) \to \Z/n\Z\big)$, and the required diagram commutes.
\end{demo}

\section{Proof of the main theorem} \label{mainsect}
Now we are going to prove our main theorem.  We change the notation for better compatibility
with the notation of section~\ref{classfieldsect}.4 and  use the letter $\X$ (instead of $X$)
for the scheme in question. We assume for the rest of this section that $\X$ is a
d-dimensional connected regular, proper and flat scheme over $\Spec(\Z)$ such that $X_\eta= \X \otimes_\Z \Q$ is projective over $\Spec(\Q)$. Let, as in section~\ref{classfieldsect}.4,
$k$ be the algebraic closure of $\Q$ in the function field of $X_\eta$ and let $S=\Spec(\Cal
O_k)$. The structural morphism $\X \to \Spec(\Z)$ factors through $S$ and $X_\eta$ is
geometrically irreducible as a variety over $k$. Let $S_f$ be the set of closed points of
$S$ and let $S_\infty$ be the set of archimedean places of the number field $k$. For $v \in
S_f$ let $Y_v = \X \otimes_S v$ be the special fibre of $\X$ over $v$. For $v \in S_\infty
\cup S_f$ let $k_v$ be the algebraic closure of $k$ in the completion of $k$ at $v$ and
$X_v=X_\eta \otimes_k k_v$.

Let $\D=\D_1 + \cdots +\D_r$ be a sum of {\em vertical} divisors on $\X$ and $\Zt=\supp(\D)$.
We choose a sufficiently small open subscheme $U\subset S$ such that $\X_U$ is disjoint to
$\Zt$. For $v \in S_f$, we denote the base change of $\X,\D,\Zt$ to $\Spec(\O_v)$ by
$\X_v,\D_v,\Zt_v$. Let $P_U$ denote the set of places of $k$ (including the archimedean ones)
which are not in $U_0$ and let $P_U^a$ and $P_U^f$ be the subset of archimedean and
nonarchimedean places in $P_U$, respectively. Recall the id\`{e}le group
(section~\ref{classfieldsect}.4)
\[
I(\X/U)= \left(\prod_{v \in P_U} \SK_1(X_v) \right) \times \left(\bigoplus_{v\in U_0}
\CH_0(Y_v) \right) .
\]
For $v \in P_U$ we endow the fields $k_v$ with the restriction of the natural topology of the
completion of $k$ at $v$ to $k_v$. This induces a natural topology on $\SK_1(X_v)$ for $v
\in P_U$ (cf. section~\ref{classfieldsect}), and (giving $\CH_0(Y_v)$ the discrete topology
for $v\in U_0$), on $I(\X/U)$. Furthermore, recall that the id\`{e}le class group $C(\X/U)$ is
defined as the quotient of $I(\X/U)$ by the image of the natural diagonal map $\SK_1(X_\eta) \to
I(\X/U)$.

\begin{defi}
We define $U^t_\D C(\X/U)$ as the image of the natural composite map
\[
\Bigg(\bigoplus_{v \in P_U^a}\; \bigoplus_{x \in (X_v)_0} k(x)^\times \;\; \oplus
\bigoplus_{ v \in P_U^f}\; U^t_\D(\X_v)\Bigg) \to I(\X/U) \to C(\X/U).
\]
\end{defi}
For $v \in P_U^f$, the group $U^t_\D(\X_v)$ contains the subgroup
$
\bigoplus_{x \in (X_v)_0} U^1(k(x)^\times)
$ (principal units with respect to the structure of $k(x)$ as a discrete valuation field).
Therefore $U^t_\D C(\X/U)$ is open in $C(\X/U)$.

\begin{theorem} \label{factorprop}
Let $\X'$ be the normalization of $\X$ in a finite cyclic extension of its function field
such that $\X'_U/\X_U$ is \'{e}tale. Then $\X'/\X$ is at most tamely ramified along $\D$ and
every real point of $\X-\Zt$ splits completely in $\X'$ if and only if the associated
character
\[
\chi_{\X'}: C(\X/U) \to \Q/\Z
\]
annihilates the subgroup $U^t_\D C(\X/U)$.
\end{theorem}

\begin{demo}{Proof:}
Firstly we assume that $\X'/\X$ is \'{e}tale over $\X-\Zt$, at most tamely ramified along $\D$ and that every real point of $\X-\Zt$ splits completely in $\X'$. The reciprocity map $\tau: C(\X/U)\to \pi_1(\X_U)^\ab$ of \cite{Sa} (cf.\ section~\ref{classfieldsect}.4) is constructed using the reciprocity maps $k(x)^\times \to \pi_1(x)^\ab$ of the fields $k(x)$ for $x \in (X_v)_0$, $v \in P_U$. The results of
section~\ref{classfieldsect}.2 show that the image of
\[
\bigoplus_{v \in P_U^a}\; \bigoplus_{x \in (X_v)_0} k(x)^\times
\]
in $C(\X/U)$ is annihilated by $\chi_{\X'}$.  For $v \in P_U^f$, we have a commutative diagram
\diagram{ccc}
\ds\bigoplus_{x \in (X_v)_0} k(x)^\times &\lang& \pi_1^t(\X_v,\D_v)^\ab\\
\mapd{}&&\mapd{}\\
C(\X/U)&\mapr{\tau}&\pi_1^t(\X,\D)^\ab
\enddiagram 
By proposition~\ref{utd}, $\chi_{\X'}$ annihilates the image of $U^t_{\D_v}(\X_v)$ in $C(\X/U)$. Therefore $\chi_{\X'}$ annihilates $U_\D^tC(\X/U)$.

\medskip
Now assume that $\chi_{\X'}$ annihilates $U^t_\D C(\X/U)$. We may assume that $\X' \to \X$ is
cyclic of prime power degree, say of order $p^\alpha$. Let $v \in P_U^f$  and $x \in
(X_v)_0$. Let $y \in Y_v$ be the unique point to which $x$ specializes. The natural map
$k(x)^\times \to U_\D^tC(\X/U)$ annihilates the subgroup  $U^1k(x)^\times$ or
$U^0k(x)^\times$ of $k(x)^\times$ if $y \in\Zt_v$ or $y\notin \Zt_v$, respectively. By local
class field theory,  the character $\chi_{\X'}$ induces a cyclic field extension of the
henselian field $k(x)$ which is tamely ramified or unramified, respectively.  In order to
show that $\X'/\X$ is \'{e}tale over $\X -\Zt$, it suffices (purity of the branch locus) to show
that it is unramified at every prime divisor $E$ not contained in $\Zt$. Since $\X'_U/\X_U$
is \'{e}tale, we may assume that $E\subset Y_v$ for some $v \in P_U^f$. The branch locus is
closed, and so it suffices to find a closed unramified point on $E$. This is easy, as for
any closed point $y\in \Zt_v$ which is a regular point of the reduced subscheme $Y_{v,
\red}$ of $Y_v$, the required property follows from lemma~\ref{saitolemma}. It remains to
show the tameness along $\Zt$. By lemma \ref{tamebyrank1}, it suffices to show that $\X'/\X$
is tamely ramified at all generic points of $\Zt$.  Tameness at points of residue
characteristic different to $p$ is immediate. Let $E$ be a prime divisor contained in $\Zt
\cap Y_v$ for an $v \in P_U$, $v |p$. We will show that $\X'/\X$ is unramified along $E$.
Let $y$ by a closed regular point of $E$. For every point $x \in (X_v)_0$ that specializes
to $y$, $k(x)$ is a henselian field of characteristic zero and of residue field
characteristic $p$. Since $\chi_{\X'}$ annihilates $U_\D^tC(\X/U)$, the associated character
on $k(x)^\times$ annihilates $U^1k(x)^\times$. It is therefore tamely ramified and of
$p$-power order, hence unramified. Using lemma~\ref{saitolemma}, we conclude that $\X'/\X$
is \'{e}tale at $y$, and therefore also at the generic point of $E$. Hence $\X'/\X$ is tame
along $\Zt$. The remaining assertion concerning the real points follows in a straightforward
manner from the results of section~\ref{classfieldsect}.2. This concludes the proof.
\end{demo}

Theorems \ref{factorprop}, \ref{endlichkeit} and the results of section A.4 imply the
\begin{corol}\label{taumodn} Putting $C^t_\D(\X/U)\stackrel{\mathrm{df}}{=} C(\X/U)/U^t_\D(\X/U)$, we obtain a surjective map
\[
\tau':C^t_\D(\X/U) \surjr{} \tilde \pi_1^t(\X,\D)^\ab.
\]
For all $n\in \N$ the map
\[
\tau'/n: C^t_\D(\X/U)/n \liso \tilde \pi_1^t(\X,\D)^\ab/n
\]
is an isomorphism of finite abelian groups.
\end{corol}

For the proof of our main theorem, we need the following approximation result:

\begin{lemma}\label{kaprox}
Assume that $\D$ is the sum of the irreducible components of the fibres over the points in $P^f_U$.  Consider the diagonal map
\[
\diag_2:  \bigoplus_{x \in (X_\eta)_1} K_2(k(x)) \lang \bigoplus_{v\in P^f_U} \bigoplus_{x \in (X_v)_1} K_2(k(x))
\]
and the boundary map 
\[
\delta_{\loc}: \ds\bigoplus_{v \in P_U^f}\; \bigoplus_{x \in (X_v)_1} K_2(k(x)) \lang
\ds\bigoplus_{v \in P_U^f}\;\ds\bigoplus_{x \in (X_v)_0} k(x)^\times \; .
\]
Let $c \in \bigoplus_{v\in P^f_U} \bigoplus_{x \in (X_v)_1} K_2(k(x))$ be given.  Then there exist an $e \in \bigoplus_{x \in (X_\eta)_1} K_2(k(x))$ with 
\[
\delta_\loc(c -\diag_2(e)) \in \bigoplus_{v \in P_U^f} \bigoplus_{y \in (Y_v)_0}
U^t_\D(y) \subset \ds\bigoplus_{v \in P_U^f}\;\ds\bigoplus_{x \in (X_v)_0} k(x)^\times.
\]
\end{lemma}

\begin{demo}{Proof:} First we choose a finite set of points $T \subset (X_\eta)_1$  such that $c_{x}=0$, for all $x \in (X_v)_1$, $v \in P^f_U$, mapping to $X_1 -T$. Then we choose local parameters $z_1,\ldots,z_r$ defining $\D_1,\ldots,\D_r$ in an affine neighborhood of $\overline{T}\cap \Zt$. 
Now let $x\in T$ and let $W$ be the closure of $\{x\}$ in $\X$. For $v \in P_U^f$ we consider those $\D_i$, $1\leq i\leq r$, which are in  $Y_v$.  The irreducible components of the pre-image of $W\cap \D_i$ in the normalization $\tilde W$ of $W$ induce discrete valuations $w_1,\ldots,w_s$ on $k(x)$. The restriction of all these valuations to $k$ is $v$. The natural maps $k(x)\otimes_k k_v \to k(x)_{w_i}$ for $i=1,\ldots, s$, induce a partition of the set $\{w_1,\ldots,w_s\}$ with respect to the \lq underlying\rq\ point in $(X_v)_1$. Applying this construction all $v \in P_U^f$  we end up with a finite number of valuations of $k(x)$ for each $x \in T$. 

\medskip
Now we put $e_x=0$ for $x \notin T$ and we choose $e_x \in K_2(k(x))$ sufficiently close to $c_{x_i}$ ($x_i \in (X_v)_1$ mapping to $x$) with respect to the valuations constructed above. By near we mean near in the sense that we approximate the entries of the symbols.  Then the image of $\diag_2(e)-c$ under the boundary map
\[
\ds\bigoplus_{v \in P_U^f}\; \bigoplus_{x \in (X_v)_1} K_2(k(x)) \lang \bigoplus_{v\in P_U^f} \bigoplus_{y \in (Y_v)_1} k(y)^\times
\]
is zero and the same holds for the images of the elements $(\diag_2(e)-c)\cup z_i$, $i=1,\ldots,r$, under the map
\[
\ds\bigoplus_{v \in P_U^f}\; \bigoplus_{x \in (X_v)_1} K_3(k(x)) \lang \bigoplus_{v\in P_U^f} \bigoplus_{y \in (Y_v)_1}K_2(k(y)).
\]
Consider the complex
\[
\ds
\bigoplus_{\renewcommand{\arraystretch}{0.8}\begin{array}{c} \sst v \in P_U^f \\
\sst x \in (X_v)_1
\end{array}
\renewcommand{\arraystretch}{1}}
 K_{j+2}(k(x)) \lang \!\!\!\ds\bigoplus_{\renewcommand{\arraystretch}{0.8}\begin{array}{c} \sst v \in P_U^f \\
\sst x \in (X_v)_0
\end{array}
\renewcommand{\arraystretch}{1}}
 K_{j+1}(k(x)) \oplus \bigoplus_{\renewcommand{\arraystretch}{0.8}\begin{array}{c} \sst v \in P_U^f \\
\sst x \in (Y_v)_1
\end{array}
\renewcommand{\arraystretch}{1}} K_{j+1}(k(y))\lang \!\!\!
\bigoplus_{\renewcommand{\arraystretch}{0.8}\begin{array}{c} \sst v \in P_U^f \\
\sst x \in (Y_v)_0
\end{array}
\renewcommand{\arraystretch}{1}}
K_{j}(k(y)).
\]
For $j=0$ we conclude that $\delta_\loc(\diag_2(e)-c)$ is in the kernel of the boundary map
\[
\ds\bigoplus_{v \in P_U^f}\; \bigoplus_{x \in (X_v)_0} k(x)^\times \lang \bigoplus_{v\in P_U^f} \bigoplus_{y \in (Y_v)_0} \Z.
\]
For $j=1$ and $1\leq i\leq r$, $\delta_\loc(\diag_2(e)-c)\cup z_i$ is in the kernel of the map
\[
\ds\bigoplus_{v \in P_U^f}\; \bigoplus_{x \in (X_v)_0} K_2(k(x)) \lang \bigoplus_{v\in P_U^f} \bigoplus_{y \in (Y_v)_0} k(y)^\times \mapr{pr}  \bigoplus_{y \in (\Zt_i)_0} k(y)^\times.
\]
By lemma~\ref{commdiag} this implies that $\delta_\loc(\diag_2(e)-c)$ is in the kernel of
\[
\phi:\ker\left(\ds\bigoplus_{v \in P_U^f}\; \bigoplus_{x \in (X_v)_0} k(x)^\times \lang \bigoplus_{v\in P_U^f} \bigoplus_{y \in (Y_v)_0} \Z\right) \lang \bigoplus_{i=1}^r \bigoplus_{y \in (\Zt_i)_0} k(y)^\times.
\]
By definition, this means that $\delta_\loc(\diag_2(e)-c)\in \bigoplus_{v \in P_U^f} \bigoplus_{y \in (Y_v)_0} U^t_\D(y)$.
\end{demo}

Now we are going to prove the following theorem (which is slightly sharper  than theorem 1 from of introduction).
\begin{theorem} \label{sharper}
Let $\X$ be a connected regular, proper and flat scheme over $\Spec(\Z)$ such that $X= \X
\otimes_\Z \Q$ is projective over $\Spec(\Q)$. Let $\D$ be a divisor on $\X$ whose 
vertical components satisfy condition $(\mu)$ of definition~\ref{condmu} and put $\Zt=\supp(\D)$.  Then the map
\[
r: \bigoplus_{x \in (\X-\Zt)_0} \Z \lang \tilde \pi_1^t (\X,\D)^\ab; \quad 1_x\mapsto \Frob_x
\]
induces an isomorphism of finite abelian groups
\[
\rec_{\X,\D}: \CH_0(\X,\D) \to \tilde \pi_1^t
(\X,\D)^\ab.
\]
\end{theorem}

\begin{demo}{Proof:} By proposition~\ref{horver1}(ii) and \ref{hotver2}, we may assume that $\D$ is vertical. We proceed in several steps.

\medskip\noindent
{\it Step 1. Existence of $\rec$}: We consider the possibly larger divisor $\D'$ which is
the sum of all vertical prime divisors lying over points in $P_U^f$. With $\Zt'=\supp(\D')$
we have $\X -\Zt'=\X_U$.  We show that the natural map
\[
r': \bigoplus_{x \in (\X-\Zt')_0} \Z \lang C^t_\D(\X/U)
\]
factors through $\CH_0^{\D'}(\X,\D)$. Then proposition~\ref{horver1} (i) and
corollary~\ref{taumodn} imply the existence of
\[
\rec: \CH_0(\X,\D)\cong \CH_0^{\D'}(\X,\D) \lang C^t_\D(\X/U) \mapr{\tau'} \tilde \pi_1^t
(\X,\D)^\ab.
\]
It remains to show that $r'$ annihilates the image of
\[
R_{\X,\D,\D'}= \ker \Big( \ker \big(\bigoplus_{x \in (\X-\Zt')_1} k(x)^\times \to
\bigoplus_{y \in \Zt'_0} \Z \big) \lang  \bigoplus_{z \in {(\Zt^+)}_0} k(z)^\times \Big)
\]
under the divisor map.  We have a decomposition
\[
(\X-\Zt')_1= (X_\eta)_0 \cup \bigcup_{v\in U_0} (Y_v)_1.
\]
For any point $x \in (Y_v)_1$, $v \in U_0$, the image of $k(x)^\times$ in
$C^t_\D(\X/U)$ is trivial by definition. We obtain an anti-commutative diagram 
\diagram{ccc}
\ds\bigoplus_{x \in (\X - \Zt')_1} k(x)^\times & \mapr{\div} & \ds\bigoplus_{x \in (\X-\Zt')_0} \Z\\
\mapd{\psi}&&\mapd{r'}\\
\ds\bigoplus_{v\in P_U^f} \bigoplus_{y\in (Y_v)_0} \big(\bigoplus_{\renewcommand{\arraystretch}{0.8}\begin{array}{c} \sst x \in (X_v)_0 \\
\sst x \to y
\end{array}
\renewcommand{\arraystretch}{1}} k(x)^\times \big)/U_\D^t(y) &\lang&C_\D^t(\X,U)
\enddiagram
in which $\psi$ is the diagonal map on $k(x)^\times$ for $x \in (X_\eta)_0$ and the zero map  for $x \in (Y_v)_1$, $v \in U_0$.  Hence $\psi$, and therefore also $r'\circ \div$ annihilates
$R_{\X,\D,\D'}$, which shows the existence of $\rec$ and, moreover, the factorization
\[
\rec: \CH_0(\X,\D) \mapr{\alpha} C^t_\D(\X/U) \mapr{\tau'} \tilde \pi_1^t
(\X,\D)^\ab.
\]

\medskip\noindent
{\it Second Step: $\tau'$ is an isomorphism of finite abelian groups}:  By
corollary~\ref{taumodn} it suffices to show that $C^t_\D(\X/U)$ is finite. The image of diagonal map
$\SK_1(X_\eta) \to \prod_{v \in P_U} \SK_1(X_v)$ is dense, and $U^t_\D C(\X/ U)$ is an open
subgroup of $C(\X/U)$. Therefore, every element in $C^t_\D(\X/U)$ can be represented by an
element in $\bigoplus_{v \in U_0} \CH_0(Y_v)$ which shows that $\alpha$ is surjective. Finally, the finiteness of $\CH_0(\X,\D)$ shows the statement of step 2.

\medskip\noindent
{\it Last step: $\alpha$ is an isomorphism:}  We show that $\alpha': \CH_0^{\D'}(\X,\D) \to C^t_\D(\X/U)$ is injective, which finishes the proof in view of proposition~\ref{horver1} (i). Consider the homomorphisms

\diagram{rrcl}
\div_{\gen}:& \ds\bigoplus_{x \in (X_\eta)_0} k(x)^\times& \lang& \ds\bigoplus_{x \in (\X-\Zt)_0} \Z\;,\\
\div_{\sp}: &\ds\bigoplus_{v \in U_0} \ds\bigoplus_{y \in (Y_v)_1} k(y)^\times& \lang
&\ds\bigoplus_{x
\in (\X-\Zt)_0} \Z\;,\\
\delta_{\loc}: &\ds\bigoplus_{v \in P_U^f}\; \bigoplus_{x \in (X_v)_1} K_2(k(x))& \lang
&\ds\bigoplus_{v \in P_U^f}\;\ds\bigoplus_{x \in (X_v)_0}
k(x)^\times\;,\\
 \delta_{\gl}:& \ds\bigoplus_{x \in (X_\eta)_1} K_2(k(x)) &\lang &\ds\bigoplus_{x \in (X_\eta)_0}
k(x)^\times\;,
\enddiagram
which are induced by the boundary maps of the respective Quillen spectral sequences.
Further, consider the diagonal maps
\diagram{rrcl}
\diag_1:& \ds\bigoplus_{x \in (X_\eta)_0} k(x)^\times&\lang \ds\bigoplus_{v \in P_U^f
}\;\bigoplus_{x \in
(X_v)_0} k(x)^\times\,,\\
 \diag_2:&\ds\bigoplus_{x \in (X_\eta)_1}K_2(k(x))&\lang
 \bigoplus_{v \in P_U^f}\; \ds\bigoplus_{x \in (X_v)_1} K_2(k(x)).
\enddiagram
Note that $\diag_1\circ\delta_{\gl}=\delta_{\loc}\circ \diag_2$. Further note that
\[
\bigoplus_{v \in P_U^f }\;\bigoplus_{x \in (X_v)_0} k(x)^\times= \bigoplus_{v \in P_U^f}\;
\bigoplus_{y \in (Y_v)_0}
\bigoplus_{\renewcommand{\arraystretch}{0.8}\begin{array}{c} \sst x \in (X_v)_0 \\
\sst x \to y
\end{array}
\renewcommand{\arraystretch}{1}}
k(x)^\times\;.
\]
and, with respect to this identity,
\[
\bigoplus_{v \in P_U^f } R_{\X_v,\D_v,\D_v'} =  \bigoplus_{v \in P_U^f } \bigoplus_{y\in (Y_v)_0}
U_\D^t(y).
\]
where $R_{\X_v,\D_v,\D_v'}$ is the group defined before proposition~\ref{horver1} for the triple $(\X_v,\D_v,\D_v')$.

\bigskip\noindent
Now let $\sum n_i P_i \in \bigoplus_{x \in (\X-\Zt')_0} \Z$ be a zero-cycle which
represents an element in $\ker(\CH_0^{\D'}(\X,\D) \to C^t_\D (\X/U))$. By definition, there exist
elements
\[
a \in \bigoplus_{x \in (X_\eta)_0} k(x)^\times, \; b \in \bigoplus_{v \in U_0} \bigoplus_{y \in
(Y_v)_1} k(y)^\times; \; c \in \bigoplus_{v \in P_U^f}\; \bigoplus_{x \in (X_v)_1} K_2(k(x)),
\]
such that the following conditions (1) and (2) hold.
\[
\sum n_i P_i = \div_{\gen}(a) + \div_{\sp}(b), \leqno{(1)}
\]
\[
\diag_1(a) -\delta_{\loc}(c) \in \bigoplus_{v \in P_U^f }\; \bigoplus_{y \in (Y_v)_0}
U^t_\D(y).\leqno{(2)}
\]
If already $\diag_1(a)$ would be in $\bigoplus_{v \in P_U^f} \bigoplus_{y \in (Y_v)_0}
U^t_\D(y)$, then  $\div_{gen}(a)$ and hence also $\div_{gen}(a)+\div_{sp}(b)$ would be a
relation in $\CH_0^{\D'}(\X,\D)$ showing that the class of $\sum \alpha_i P_i$ is trivial in
$\CH_0^{\D'}(\X,\D)$. We will achieve this by replacing $a$ by $a-\delta_{\gl}(e)$, where $e$
is a suitable chosen element in\break $\bigoplus_{x \in (X_\eta)_1}K_2(k(x))$. In view of the complex
\[
\bigoplus_{x \in \X^{d-2}} K_2(k(x)) \lang \bigoplus_{x \in \X^{d-1}} k(x)^\times \lang
\bigoplus_{x \in \X^{d}} \Z
\]
this does not change $\div_{\gen}(a)$ modulo the image of $\div_{\sp}$.
Now we construct $e$. By lemma~\ref{kaprox} applied to $\D'$, we find $e$ such that 
\[
\delta_\loc(\diag_2(e)-c)\in \bigoplus_{v \in P_U^f} \bigoplus_{y \in (Y_v)_0} U^t_{\D'}(y) \subset \bigoplus_{v \in P_U^f} \bigoplus_{y \in (Y_v)_0} U^t_\D(y) 
\]
Setting, $a'=a-\delta_{\gl}( e)$, we obtain
\[
\begin{array}{rcl}
\diag_1(a')&=& \diag_1(a)-\diag_1(\delta_{\gl}(e))\\
&=& \diag_1(a)-\delta_{\loc}(c)+\delta_{\loc}(c-\diag_2(e)) \in \bigoplus_{v \in P_U^f }
R_{\X_v,\D_v,\D_v'}.
\end{array}
\]
Therefore $a' \in R_{\X,\D,\D'}$ and we conclude that the class of $\sum \alpha_i P_i$ is
zero in $\CH_0^{\D'}(\X,\D)$.
\end{demo}

\vskip.8cm The following variant of theorem~\ref{sharper} describes the full
group $\pi_1^t(\X,\D)^\ab$ and not its quotient $\tilde\pi_1^t(\X,\D)^\ab$. Of course, this
yields nothing new if $\X$ has no real points. We define a modified relative Chow group by
allowing only such relations $\div (f)$ of functions $f$ which are positive at every
real-valued point. More precisely, for a point $x \in \X_1$ we put
\[
k(x)_+^\times=\{f \in k(x)^\times\,|\, \iota(f)>0 \hbox{ for every embedding } \iota: k(x)
\to \R \}.
\]
If the global field $k(x)$ is of positive characteristic or a totally imaginary number
field, then $k(x)^\times_+ =k(x)^\times$.

\begin{defi}
The group $\widetilde \CH_0(\X,\D)$ is the quotient of $\ds\bigoplus_{x\in (\X-\Zt)_0}\Z$ by the
image of
\[
\tilde R_{\X,\D}=\ker \Big( \ker \big(\bigoplus_{x \in (\X-\Zt)_1} k(x)^\times_+ \to \bigoplus_{y
\in \Zt_0} \Z \big) \lang  \bigoplus_{z \in \Zt^+_0} k(z)^\times \Big)
\]
under the divisor map.
\end{defi}

A slight modification of the proof of theorem~\ref{sharper}, which we leave to the reader, shows the

\begin{theorem} Let $\X$ be a regular connected scheme, flat and proper over $\Spec(\Z)$ such that its generic fibre $\X \otimes_\Z \Q$ is projective over $\Q$. Let $\D$ be a divisor on $\X$ whose vertical components satisfy condition $(\mu)$. Then there exists a natural isomorphism of finite abelian groups
\[
\widetilde \rec_{\X,\D}: \widetilde \CH_0(\X,\D) \liso \pi_1^t(\X,\D)^\ab\;.
\]
\end{theorem}

\appendix
\section{Appendix: Results from class field theory}\label{classfieldsect}

For reference, to fix notation and to deduce some  corollaries,  we recall several results
from the higher dimensional local and global class field theory (mainly due to the work of K.
Kato and S. Saito). Let $X$ be a connected noetherian scheme. For $j \in \N$, let
\[
X_j=\left\{x \in X\,|\, \dim \left(\overline{\{x \}}\right)=j \right\}, \quad X^j=\left\{x
\in X\,|\, \codim_X \left(\overline{\{x \}}\right)=j \right\}.
\]
For an integer $i\geq 0$, the localization theory in Quillen $K$-theory on $X$ gives rise to
a homomorphism
\begin{equation}\label{divmap}
 \partial_i: \bigoplus_{y \in X_1} K_{i+1}(k(y)) \lang \bigoplus_{x \in X_0} K_i(k(x)).
\end{equation}
We define $\SK_i(X)$ to be the cokernel of this map $\partial_i$. For $i=0,1$, theorems 2.2
and 2.6 of \cite{Sn} imply the following explicit descriptions of $\partial_i$: The
restriction of $\partial_i$ to the $y$-component $K_{i+1}(k(y))$ is the composite of the
following homomorphisms:
\[
K_{i+1}(k(y)) \mapr{\partial_i'} \bigoplus_{x' \in \tilde{Y}_0} K_i(k(x')) \mapr{N} \bigoplus_{x
\in Y_0} K_i(k(x)) \lang \bigoplus_{x\in X_0} K_i(k(x))
\]
where $Y$ is the closure of $\{y\}$ in $X$, $\tilde{Y}$ is the normalization of $Y$ and the map
$\partial_i'$ is defined as
\[
\partial_i'=\left\{
\begin{array}{l}
\ds\sum_{x'\in \tilde{Y}_0} \ord_{x'} \hbox{ \rm for } i=0,\\
\ds \sum_{x'\in \tilde{Y}_0} d_{x'} \hbox{ \rm for } i=1,
\end{array}
 \right.
\]
where $d_{x'}$ is defined for $t=(f,g)\in K_2(k(y))$ by
\[
d_{x'}(t)=(-1)^{\ord_{x'}(f)\ord_{x'}(g)} f^{\ord_{x'}(g)}g^{-\ord_{x'}(f)}.
\]
Finally, the map $N$ is defined as follows: For $x'\in \tilde{Y}_0$ and $x\in Y_0$, the
$(x',x)$-component of $N$: $K_i(k(x'))\to K_i(k(x))$ is the zero map if $x'$ does not lie
over $x$ ; otherwise it is multiplication by $[k(x'):k(x)]$: $\Z \to \Z$ if $i=0$, and the
usual norm homomorphism $k(x')^\times \to k(x)^\times$ if $i=1$.

\bigskip\noindent
Assume that $X$ is a variety over a field $k$ and that $k$ (and hence also every finite
extension of $k$) carries a natural topology. Then we give $\SK_0(X)$ the discrete topology
and we endow $\SK_1(X)$ with the finest topology such that for every $x\in X_0$ the natural
homomorphism
\[
k(x)^\times \lang \SK_1(X)
\]
is continuous.

\bigskip\noindent
\underline{A 1. The reciprocity homomorphism for an arithmetic scheme.}

\medskip\noindent

Let $X$ be a connected scheme of finite type over $\Spec(\Z)$. For each closed point $x \in
X_0$, $k(x)$ is a finite field, so there is an isomorphism
\[
\hat \Z \liso \pi_1^\ab(x)=Gal(k(x)^{sep}|k(k)),
\]
which sends $1\in \hat \Z$ to the Frobenius $f_x$ over $k(x)$. We define the Frobenius
element $F_x$ of $x$ in $\pi_1^\ab(X)$ to be the image of $f_x$ under the natural
homomorphism $\pi_1^\ab(x) \to \pi_1^\ab(X)$. The assignment $1 \mapsto F_x$, $x \in X_0$,
defines a homomorphism

\begin{equation} \label{recmap}
\bigoplus_{x \in X_0} \Z \lang \pi_1^\ab (X).
\end{equation}

By the results of Lang \cite{La}, we have the
\begin{lemma}\label{dense}
If $X$ is irreducible and the reduced subscheme $X_{red}$ of $X$ is normal, then the map
\eqref{recmap} has a dense image.
\end{lemma}

Now we consider the abelianized modified fundamental group $\tilde\pi_1^\ab(X)$ which is the
unique quotient of $\pi_1^\ab(X)$ classifying finite abelian \'{e}tale coverings of $X$ in which
every real-valued point of $X$ splits completely.

\begin{lemma}{\rm (\cite{Sa}, lemma 2.4)}
If $X$ is proper over $\Spec(\Z)$, then the composite of the map \eqref{recmap} and the
natural surjection $\pi_1^\ab (X) \to \tilde \pi_1^\ab (X)$ annihilates the image of
\[
\partial: \bigoplus_{x \in X_1} k(x)^\times \lang \bigoplus_{x \in X_0} \Z,
\]
which is the map \eqref{divmap} for $i=0$.
\end{lemma}

Consequently, we obtain a natural map
\[
\SK_0(X) \lang \tilde \pi_1^\ab(X),
\]
which is called the reciprocity map for $X$.

\bigskip\noindent
\underline{A 2. Class field theory for varieties over $\R$ and $\C$.}

\medskip\noindent
We follow \cite{Sa}, \S4. Let $k^*$ be $\R$ or $\C$ and let $X^*$ be a connected, proper and
smooth scheme over $k^*$. For $x \in X^*_0$, $k(x)\cong \R$ or $\C$ and we have a canonical
surjection
\[
k(x)^\times \lang \pi_1^\ab(x)=Gal({k(x)}^{sep}|k(x)).
\]
The maps $\pi_1^\ab(x) \to \pi_1^\ab(X^*)$ for $x \in X^*_0$ induce a well-defined continuous
homomorphism
\begin{equation}\label{eq1}
\tau: \SK_1(X^*) \lang \pi_1^\ab(X^*),
\end{equation}
which is the zero map if $k^*=\C$ or $X^*(\R)=\varnothing$.  Let $k\subset k^*$ be a
subfield and suppose that the following conditions (a)--(c) are satisfied

\begin{itemize}
 \item[(a)] $k$ is dense in $k^*$ for the usual topology of $k^*$.
 \item[(b)] If $k^*=\C$, $k$ is algebraically closed.
 If $k^*=\R$, the algebraic closure $\bar k$
of $k$ is a quadratic extension of $k$.
 \item[(c)] There is a proper, smooth scheme $X$ over $k$ such that $X \otimes_k k^* \cong X^*$.
\end{itemize}

\noindent By the proper base change theorem, the natural map $\pi_1^\ab(X^*)\to \pi_1^\ab(X)$
is an isomorphism and the image of the natural homomorphism
\begin{equation}\label{eq2}
\SK_1(X) \lang \pi_1^\ab(X)
\end{equation}
coincides with that of the map (\ref{eq1}). Slightly more general, let $V \subset X$ be a
nonempty open subscheme and $V^*=V\otimes_k k^*$. Using the smooth base change theorem
instead (we are in characteristic zero), we have again the isomorphism $\pi_1^\ab(V^*) \liso
\pi_1^\ab(V)$ and, sending a point $x \in V(\R)=V^*(\R)$ to the image of the unique
nontrivial element of $\pi_1^\ab(x)\cong Gal(\C|\R)$ in $\pi_1^\ab(V)$, we obtain a natural
map
\[
i: V(\R) \lang \pi_1^\ab(V).
\]
The map $i$ is locally constant on $V(\R)$ and has finite image  (see \cite{Sa}, lemma 4.8
for the case $V=X$, the proof in the general case is the same). If $\chi\in
H^1_\et(V,\Q/\Z)$ corresponds to a cyclic \'{e}tale covering $ V' \to V$, then a point $x
\in V(\R)$ splits completely in $V'$ if and only if $\chi$ is trivial on $i(x)$. In
particular, the subset of points in $V(\R)$ which split completely in $V'$ is (norm)
open and closed in $V(\R)$.

\bigskip\noindent
\underline{A 3. Class field theory of schemes over henselian discrete valuation fields.}

\medskip\noindent
We follow \cite{Sa}, \S3. Let $\O_k$ be a henselian discrete valuation ring with finite
residue field $F$ and quotient field $k$. Let $X$ be a proper smooth scheme over $k$. For $x
\in X_0$, $k(x)$ is a finite extension of $k$, so that the local class field theory for
$k(x)$ gives us a canonical homomorphism
\[
\rho_x : k(x)^\times \lang \pi_1^\ab(x) \lang \pi_1^\ab (X).
\]
The sum of these $\rho_x$ is trivial on the diagonal image of $\bigoplus_{x \in
X_1}K_2(k(X))$ and we obtain a canonical morphism
\begin{equation}\label{eq3}
 \tau: \SK_1(X) \lang \pi_1^\ab(X).
\end{equation}
Giving each $k(x)$, $x\in X_0$ the usual adic topology, we obtain a natural topology on
$\SK_1(X)$ (see the beginning of this section). The map  $\tau$ is continuous with respect to
this topology and the natural profinite topology on $\pi_1^\ab(X)$. For every natural number
$n$ which is prime to char$(k)$, the subgroup $n\cdot\SK_1(X)$ is open in $\SK_1(X)$. In
particular, if char$(k)$=0, then every homomorphism $\SK_1(X) \to \Q/\Z$ with finite image is
continuous.

\medskip
Now assume that $X$ has a proper regular model $\X$ over $\O_k$ and let $Y$ be its special
fibre. By \cite{Sa}, lemma 3.11, we have a natural continuous and surjective homomorphism
$\partial: \SK_1(X) \to \SK_0(Y)$ and a commutative diagram
 \diagram{ccc}
 \SK_1(X)&\mapr{\tau} & \pi_1^\ab(X)\\
 \mapd{\partial}&&\mapd{sp}\\
 \SK_0(Y) &\mapr{\phi} &\pi_1^\ab (Y),
 \enddiagram
where $sp$ is the specialization homomorphism (cf.\ \cite{SGA1}, X) and $\phi$ is the
reciprocity map for $Y$ as defined in paragraph 1 above. Important for us is the

\begin{prop} {\rm (\cite{Sa}, prop. 3.3)}
Let $\chi \in H^1_\et(X,\Q/\Z)$ and let $\tilde \chi: \SK_1(X) \to \Q/\Z$ be the induced
homomorphism via \eqref{eq3}. Then $\chi$ comes from the subgroup $H^1_\et(Y,\Q/\Z)\cong
H^1_\et(\X,\Q/\Z) \subset H^1_\et(X,\Q/\Z)$ if and only if $\tilde \chi$ factors though the
map $\partial$.
\end{prop}

The crucial step in the proof of the last proposition is the following lemma, which we also
need in the proof of theorem~\ref{factorprop}.

\begin{lemma}{\rm (\cite{Sa}, lemma 3.15)} \label{saitolemma}
Let $A$ be a henselian regular local ring of dimension $\geq 2$ with perfect residue field
$F$ and quotient field $K$. If\/ $\hbox{\rm char}(K)=0$ and $\hbox{\rm char}(F)=p>0$ assume
that there exists exactly one height-one prime ideal $\p \subset A$ dividing $p$ and let
$T\in A$ be an element with $\p=(T)$. In the remaining cases let $T\neq 0$ be any non-unit
of $A$. Put $U=\Spec(A[1/T])$. Then, if $\chi \in H^1_\et(U,\Q/\Z)$ induces an unramified
character $\chi_u \in H_\et^1(u,\Q/\Z)$ for each $u \in U_0$, $\chi$ comes from
$H^1_\et(\Spec(A),\Q/\Z)$.
\end{lemma}

\bigskip\noindent
\underline{A 4. Unramified class field theory of arithmetic schemes.}

\medskip\noindent
Let $\X$ be a connected regular, proper and flat scheme over $\Spec(\Z)$, and suppose that
$X_\eta= \X \otimes_\Z \Q$ is projective over $\Spec(\Q)$.  Let $k$ be the algebraic closure of $\Q$ in the function field of $X$ and put $S=\Spec(\Cal O_k)$. The structural morphism $\X
\to \Spec(\Z)$ factors through $S$ and $X_\eta$ is geometrically irreducible as a variety over $k$. Let $S_f$ be the set of closed points of $S$ and let $S_\infty$ be the set of
archimedean places of the number field $k$. For $v \in S_f$ let $Y_v = \X \otimes_S v$ be
the special fibre of $\X$ over $v$. For $v \in S_\infty \cup S_f$, let $k_v$ be the algebraic
closure of $k$ in the completion of $k$ at $v$ and $X_v=X_\eta \otimes_k k_v$.

\ms  For $v \in S_f$ we recall the continuous and surjective homomorphism \eqref{eq3} from
paragraph 3:
\[
\tau_v:\, \SK_1(X_v) \lang \SK_0(Y_v)=\CH_0(Y_v).
\]
For a nonempty open subscheme $U \subset S$ we consider the topological group
\[
I(\X/U)= \left(\prod_{v \in P_U} \SK_1(X_v) \right) \times \left(\bigoplus_{v\in U_0}
\CH_0(Y_v) \right) ,
\]
where $P_U$ denote the set of places of $k$ (including the archimedean ones) which are not in
$U_0$. There exists a natural homomorphism (\cite{Sa}, 5.3)
\[
I(\X/U) \lang \pi_1(\X_U)^{\ab},
\]
which annihilates the image of the diagonal map $\SK_1(X_\eta) \lang I(\X/U)$. It therefore
induces a continuous homomorphism (\cite{Sa}, 5.5)
\[
\tau: C(\X/U):=\coker\big(\SK_1(X_\eta) \to I(\X)\big) \lang \pi_1(\X)^\ab.
\]
By \cite{Sa}, 5.6., every subgroup of finite index in $C(\X/U)$ is open and, for every
positive integer $n$, $\tau$ induces an isomorphism of finite abelian groups
\[
C(\X/U)/n \mapr{\sim} \pi_1(\X_U)^\ab /n.
\]

\noindent \footnotesize{Alexander Schmidt, NWF I - Mathematik, Universit\"{a}t Regensburg, D-93040
Regensburg, Deutschland. email: alexander.schmidt@mathematik.uni-regensburg.de}
\end{document}